\documentclass{article}
\usepackage{amssymb}
\usepackage{amsfonts}
\usepackage{amsmath}
\usepackage{url}
\usepackage{latexsym}

\def\Gl{\mathop{\rm Gl}\nolimits}
\def\deg{\mathop{\rm deg }\nolimits}
\def\rank{\mathop{\rm rank}\nolimits}
\def\lcm{\mathop{\rm lcm }\nolimits}
\def\rev{\mathop{\rm rev }\nolimits}
\newcommand{\se}{\ensuremath{\stackrel{s.e.}{\sim}}}

\newcommand{\efe}{\mathbb F}
\newcommand{\FF}{\mathbb F}

\newcommand{\ba}{\mathbf a}
\newcommand{\bd}{\mathbf d}
\newcommand{\bg}{\mathbf g}
\newcommand{\bc}{\mathbf c}
\newcommand{\bb}{\mathbf b}
\newcommand{\bu}{\mathbf u}
\newcommand{\bv}{\mathbf v}

\newtheorem{theorem}{Theorem}[section]
\newtheorem{corollary}[theorem]{Corollary}
\newtheorem{lemma}[theorem]{Lemma}
\newtheorem{definition}[theorem]{Definition}
\newtheorem{problem}[theorem]{Problem}
\newtheorem{proposition}[theorem]{Proposition}
\newtheorem{rem}[theorem]{Remark}

\newenvironment{remark}{\begin{rem} \em}{\end{rem}}

\title{Row or column completion of polynomial matrices 
of given degree
\thanks{
This work was 
supported by grant PID2021-124827NB-I00 funded by MCIN/AEI/ 10.13039/501100011033 and by ``ERDF A way of making Europe'' by the ``European Union''.
The first and third authors were also supported 
by grant GIU21/020 funded by UPV/EHU.
}}

\author{Agurtzane Amparan\thanks{Departamento de Matem\'aticas, Universidad del Pa\'is Vasco UPV/EHU, Bilbao, Spain, {agurtzane.amparan@ehu.eus}, {silvia.marcaida@ehu.eus}}
\and Itziar Baraga\~na\thanks{Departamento de Ciencia de la Computaci\'on e I.A., 
Universidad del Pa\'{\i}s Vasco, UPV/EHU, Donostia-San Sebasti\'an, 
Spain, {itziar.baragana@ehu.eus}}
\and Silvia Marcaida\footnotemark[2]
\and Alicia Roca\thanks{Departamento de Matem\'atica Aplicada, IMM, Universitat Polit\`ecnica de Val\`encia, 46022 Valencia, Spain,   {aroca@mat.upv.es}}}

\date{}

\DeclareMathOperator{\diag}{diag}

\begin{document}

\maketitle

\begin{abstract}

We solve the problem of characterizing the existence of a polynomial matrix of fixed degree when its eigenstructure  (or part of it)    and some of its rows (columns) are prescribed.  More specifically, we present  a solution to  the row (column) completion problem of a polynomial matrix of given degree under  different   prescribed invariants: the whole  eigenstructure,  all of it but the row (column) minimal indices, and the finite and/or infinite structures.
Moreover, we characterize the existence of a polynomial matrix with prescribed  degree and eigenstructure over an arbitrary field.
\end{abstract}


{\bf AMS:}
15A18 ,15A54, 15A83, 93B18

\section{Introduction}\label{secintroduction}

 Polynomial matrices are very often used to study the dynamical behavior of systems of differential or difference equations \cite{Bhms11, Fo75, GLR85, Kail80, Rose70, TiMe01, Vard91}. These equations arise in many different scientific areas such as engineering, physics, economics and biology.  When the polynomial matrix associated with the system is regular, a closed formula for its solution can be given in terms of the finite and infinite elementary divisors of the matrix; when it is singular, the solution also depends on the left and right null spaces of the matrix, which are related to its minimal indices.

On the other hand, a very important problem in applications is the matrix completion problem. 
A matrix completion problem consists in characterizing  all or some of the invariants of a matrix with respect to a given equivalence relation, when some entries  of the matrix have been prescribed. Different types of matrices, different equivalence relations, or different   positions of the prescribed components lead to different type of problems.
For example, a variety of pencil completion problems  appear in the design of linear control systems to modify the structure of the system (see, for instance, \cite{BoDo94, DoSt19, LoMoZaZaLAA98} and the references therein.

The research on  matrix completion problems  has been very active for a long time. Since the 1950s numerous  problems in the field have been studied for constant  and polynomial matrices.
 The  literature in the area is vast, therefore the references given in this section do not intend at all  to be exhaustive, but a few examples shall be mentioned. Results have been obtained for the problem of obtaining a square constant matrix when a submatrix and similarity invariants  are prescribed \cite{CaSi92, Di74, Ol69, FaLe59, GoKaSch83, HeLMA83,  Sa79, Mi58,   Si87, Th79, ZaLAA87}, and also for rectangular  constant matrices when  invariants for the feedback equivalence relation and some of its rows or its columns are prescribed \cite{BaZa90, Do05, ZaLAA88}.
 
Some results concerning matrix pencil completion problems for the strict equivalence were obtained   in  \cite{Ba89, CaSi91,  FuSi99}. The general matrix pencil completion problem was finally stated in \cite{LoMoZaZaLAA98}, and since then enormous progresses have been made,  most of them by the same authors (see \cite{Do08, Do10,   Do13, Do22,  DoSt09, DoSt19}).

 Especially related to our work  are the references
\cite{ Do05, Do13, Sa79, Th79, ZaLAA87}.
 A fundamental result for polynomial matrices was obtained in   \cite{Sa79, Th79}. The result is known as the {\em interlacing inequalities} for invariant factors and underlies many of the above mentioned papers.
It characterizes the invariants for the unimodular equivalence of a polynomial matrix (i.e., the invariant factors)  when a submatrix is prescribed and there is no restriction on the degree of the completion. The authors also characterized the  invariants for the similarity of square constant matrices when a principal submatrix is prescribed.
The other three references are devoted to  row (column)  completion problems.
In  \cite{ZaLAA87}, the problem of prescribing the similarity invariants of a square constant matrix  when some of its rows are prescribed, is solved.
A solution to the analogous problem for rectangular constant matrices and  feedback invariants of pairs of matrices was given in   \cite{Do05}.
For matrix pencils  and prescribed invariants for the strict equivalence, a solution was obtained  in  \cite{Do08}, where  an implicit solution was provided.
In \cite{Do10} the same problem was solved in terms of a (complex) explicit solution. Later, in \cite{Do13}, the solution was simplified (see also \cite{DoSt19}).

The invariants for the strict equivalence of matrix pencils are the invariant factors  (finite structure), the infinite elementary divisors (infinite structure) and the row and column  minimal indices; i.e., the left and right minimal indices as polynomial matrices (singular structure). For polynomial matrices of arbitrary degree, these 4 types of invariants are known as the {\em eigenstructure} of the matrix. The eigenstructure of a polynomial matrix is invariant for the strict equivalence, but it does not form a complete system of invariants for such relation (see \cite{DeDoMa14}).
For  polynomial matrices the finite and infinite structures can be given in terms of the homogeneous invariant factors.

We  are interested in the following problem:
 \begin{problem}\label{problem}
Let $P(s)\in\efe[s]^{m\times n}$ be a polynomial matrix of $\deg(P(s))=d$. Find necessary and sufficient conditions for the existence of a polynomial matrix  $W(s)\in \efe[s]^{z\times n}$ of $\deg(W(s))\leq d$
such that  $\begin{bmatrix}P(s)\\W(s)\end{bmatrix}$  has prescribed eigenstructure, or part 
 of it.
\end{problem}
Since the eigenstructure consists of four types of invariants, depending on the invariant(s) prescribed we get different problems;  Problem \ref{problem} leads  actually to 15 problems.
In this paper we solve Problem \ref{problem} when the whole eigenstructure (homogeneous invariant factors and row and column minimal indices) is prescribed, therefore generalizing the result of \cite{Do13}. To achieve it, we turn a row completion problem of  a polynomial matrix into  a row completion problem of a matrix pencil, using a suitable linearization (the first Frobenious companion form) of the polynomial matrix. This allows us to use the result on row completion of  pencils  in \cite{Do13, DoSt19}. We also solve some of the particular cases appearing when part of the eigenstructure is prescribed.

Results are obtained for row completion problems, but it is a matter of transposition to obtain solutions to the corresponding column completion problems.

When no   row is prescribed, i.e., when we want to characterize the existence of a poynomial matrix with prescribed   eigenstructure, the problem was solved in \cite{DeDoVa15} for infinite fields. Here we remove this restriction and provide a solution to the problem for arbitrary fields. 

The paper is organized as follows. 
Section \ref{secpreliminaries} is a preliminary section. It is structured in three subsections. Subsection \ref{notation} contains notation and some basic results. In Subsection \ref{secpencilcompletion} we present a known result of row completion  of matrix pencils with the corresponding notation and definitions.
Subsection \ref{Frobenius} is devoted to the Frobenius companion form of a polynomial matrix.  
In Section \ref{secpolynomialprescreig} 
we characterize    the existence of a polynomial matrix  with prescribed eigenstructure in an arbitrary field.
In Section \ref{secpolcompl}  we  prove that  the row completion problem of a polynomial matrix is equivalent to the row completion problem of a pencil, via the Frobenius companion form. In different subsections we present the solution to  Problem \ref{problem} when prescribing the whole eigenstructure  (Subsection \ref{subsec_eigen}), 
the whole eigenstructure but the row (column) minimal indices  (Subsection \ref{subsec_fininfcol} (\ref{subsec_fininfrow})), 
and the finite and/or infinite elementary divisors (Subsection \ref{subsec_fininf}). 
We also include a section of conclusions and future work (Section \ref{secconclusions}).

\section{Preliminaries}
\label{secpreliminaries}

\subsection{Notation and previous results} \label{notation}
Let $\FF$ be a field and $\overline{\FF}$ its algebraic closure. $\FF[s]$ denotes the ring of polynomials in the indeterminate $s$ with coefficients in $\FF$,  $\FF(s)$ is the field of fractions of
$\FF[s]$, i.e., the field of rational functions over $\FF$, and  $\FF[s, t]$ denotes  the ring of polynomials in two
variables $s, t$ with coefficients in $\FF$. 
A polynomial in $\FF[s]$ is \textit{monic} if its leading coefficient is 1. We will say that a polynomial in $\FF[s,t]$ is {\em monic}
    if it is monic with respect to the variable $s$.
  Given two polynomials $p, q$, by $p\mid q$ we mean that $p$ is a divisor of $q$. By $\lcm(p,q)$ we mean the monic
    least common multiple of $p$ and $q$. 
   
We denote by  $\FF^{m\times n}$, $\FF[s]^{m\times n}$  and $\FF(s)^{m \times n}$  the vector spaces over $\FF$  of $m \times n$ matrices with elements in $\FF$, $\FF[s]$ and $\FF(s)$, respectively.
$\Gl_n(\FF)$ will be the general linear group of invertible matrices
in $\FF^{n \times n}$.

Let  $P(s)\in \FF[s]^{m\times n}$ be a polynomial matrix. It can be written as $P(s)=P_ds^d+P_{d-1}s^{d-1}+\cdots+P_1s+P_0\in\efe[s]^{m\times n}$, $P_i\in  \FF^{m\times n}$,  $0\leq i \leq  d$, with $P_d\neq 0$  for some integer $d$. Then, $d$ is  the {\em degree} of $P(s)$ (denoted by $\deg(P(s))$). 

The {\em normal rank} of $P(s)$, denoted by $\rank (P(s))$,  is the order of the largest non identically zero minor of $P(s)$, i.e., it is the rank of $P(s)$ considered as a matrix on  $\FF(s)$. We will refer to it as the {\em rank} of $P(s)$.

A polynomial matrix $P(s)\in \FF[s]^{m\times n}$
 is  {\em regular} if $m=n$ and $\det(P(s))$ is  not identically zero. Otherwise it is  {\em singular}.

A polynomial matrix $P(s)\in \FF[s]^{n\times n}$
 is  {\em unimodular} if $\det(P(s))$ is a non zero constant. Given a polynomial matrix  $P(s)\in \FF[s]^{m\times n}$, $r=\rank(P(s))$, there exist unimodular matrices  $U(s)\in \FF[s]^{m\times m}$ and $V(s)\in \FF[s]^{n\times n}$ such that 
$$U(s)P(s)V(s)=D(s)=\left[\begin{array}{ccc|c}
\alpha_1(s)  & & & \\
& \ddots & & 0 \\
& &  \alpha_r(s)& \\
\hline
& 0 & & 0
\end{array}\right], $$
where $\alpha_1(s)\mid  \dots \mid \alpha_r(s)$ are monic polynomials  and are called the {\em invariant factors} of $P(s)$. The matrix $D(s)$ is  the {\em Smith form} of  $P(s)$ (see, for instance, \cite{Gant66, GLR85}).

It is said that $\lambda \in \overline{\FF}\cup \{\infty\}$ is an {\em eigenvalue} of $P(s)$ if $\rank(P(\lambda))< \rank(P(s))$; here we understand that $P(\infty)=P_d$. The set of eigenvalues of $P(s)$ is the {\em spectrum} of the matrix, and we denote it by  $\Lambda(P(s))$.

Given $P(s)=P_ds^d+P_{d-1}s^{d-1}+\cdots+P_1s+P_0\in\efe[s]^{m\times n}$  with $P_d\neq 0$,  we define the {\em reversal matrix polynomial} $\rev P(t)$ as
\begin{equation} \label{reversal}
\rev P(t)=t^dP(1/t)=P_d+P_{d-1}t+\cdots+P_1t^{d-1}+P_0t^d.
\end{equation}
This matrix and $P(s)$ have the same rank but $\rev P(t)$ may be of degree smaller than $d$. It turns out that $\infty$ is an eigenvalue of $P(s)$
  if and only if $0$ is an eigenvalue of $\rev P(t)$.

Notice that the finite eigenvalues of $P(s)$ are the roots in $\overline{\FF}$ of the polynomial $\alpha_r(s)$ in the Smith form of $P(s)$.
Factorizing the invariant factors as products of irreducible polynomials over $\overline{\FF}$,  we can write
$$
\alpha_{i}(s)=\prod_{\lambda\in \Lambda(P(s))\setminus\{\infty\}}(s-\lambda)^{n_i(\lambda, P(s))}, \quad 1\leq i \leq r.
$$
The factors  $(s-\lambda)^{n_i(\lambda, P(s))}$ with $n_i(\lambda, P(s))>0$ are the {\em  elementary divisors}  of $P(s)$ over $\overline{\FF}$ corresponding to $\lambda$, and the integers
$
n_1(\lambda, P(s))\leq \dots \leq n_r(\lambda, P(s))
$
are the {\em partial multiplicities} of $\lambda$ in $P(s)$. 

The  {\em infinite elementary divisors}  of $P(s)$ and the {\em partial multiplicities of $\infty$} in $P(s)$ are the elementary divisors  of $\rev P (t)$ corresponding to $0$ and the partial multiplicities of $0$ in $\rev P(t)$, respectively. 
For simplicity, we will denote $e_i=n_i(0, \rev P(t))$ for $1\leq i \leq r$, where if
$0 \not \in \Lambda(\rev P(t))$ (i.e., $\infty  \not \in \Lambda(P(s))$),   we take $e_1=\dots=e_r=0$.

We recall now the singular structure of a polynomial matrix. Denote by $\mathcal{N}_\ell (P(s))$ and $\mathcal{N}_r (P(s))$ the {\em left} and
{\em right null-spaces} over $\FF(s)$ of $P(s)$, respectively, i.e.,
if $P(s)\in\FF[s]^{m\times n}$,
\[
\begin{array}{l}
\mathcal{N}_\ell (P(s))=\{x(s)\in\FF(s)^{m \times 1}: x(s)^TP(s)=0\},\\
\mathcal{N}_r (P(s))=\{x(s)\in\FF(s)^{n \times 1}: P(s)x(s)=0\}.
\end{array}
\]
These sets are vector subspaces
of $\FF(s)^{m \times 1}$ and $\FF(s)^{n \times 1}$, respectively. For a subspace $\mathcal{V}$ of
$\FF(s)^{m \times 1}$ it is possible to find a basis consisting of vector
polynomials; it is enough to  take an arbitrary basis and multiply each vector
by a  least common multiple of the denominators of its  entries.
The {\em order} of a polynomial basis is defined as the sum of
the degrees of its vectors  (see \cite{Fo75}). A {\em minimal basis} of  $\mathcal{V}$
is  a polynomial basis  with least order among all
polynomial bases of $\mathcal{V}$. 
The increasing ordered list of degrees of the vector polynomials  of a minimal basis   is  always the same (see \cite{Fo75}).  These degrees are called the {\em minimal indices} of $\mathcal{V}$.

 A   {\em right  (left) minimal basis} of a polynomial matrix $P(s)$ is a minimal basis of $\mathcal{N}_r (P(s))$ ($\mathcal{N}_\ell (P(s))$).
The {\em right  (left) minimal indices} of $P(s)$ are the minimal indices of
$\mathcal{N}_r (P(s))$ ($\mathcal{N}_\ell (P(s))$).  
From now on in this paper, we will work with the right (left) minimal indices ordered decreasingly, and we will refer to them as the {\em column (row) minimal indices} of $P(s)$.  Notice that a polynomial matrix $P(s)\in \FF[s]^{m\times n}$ of $\rank(P(s))=r$ has $m-r$ row  and $n-r$ column minimal indices.
The  invariant factors,  the infinite elementary divisors and the column and row minimal indices form the {\em eigenstructure} of a polynomial matrix $P(s)$ (see, for instance, \cite{DeDoVa15}).

A {\em matrix pencil} is a  polynomial matrix $P(s)\in \FF[s]^{m\times n}$ of degree at most 1. 
Two matrix pencils $A(s)=sA_1+A_0, B(s)=sB_1+B_0\in \FF[s]^{m\times n}$ are {\em strictly equivalent} ($A(s)\se B(s)$) 
if there exist invertible matrices $Q\in \Gl_m(\FF)$,   $R\in \Gl_n(\FF)$ such that 
$A(s)=QB(s)R.$
A complete system of invariants for  the strict equivalence  of matrix pencils is formed by 
the invariant factors, the infinite elementary divisors (or, equivalently, the partial multiplicities of $\infty$), and the column and row minimal indices. 
They are known as the {\em Kronecker invariants} of the pencil. 
The associated canonical form is the Kronecker canonical form. For details see  \cite[Ch. 2]{Fri16} or \cite[Ch. 12]{Gant66} for infinite fields, and  \cite[Ch. 2]{Ro03} for arbitrary fields.

\begin{remark}\label{remiedpencil}
  In the above references, and as far as we know,  the definition of the reversal  of a constant pencil is different from the definition given in (\ref{reversal})
  for polynomial matrices, which leads to  different infinite structure for constant matrices.  More specifically,
 given a pencil $A(s)=sA_1+A_0$,  the reversal pencil is generally defined as $\rev A(t)=A_1+tA_0$. According to this definition, the reversal of a constant pencil  $A_0$ is $tA_0$, therefore the partial multiplicities of $\infty$ in a constant pencil of rank  $r$ are $e_1=\dots=e_r=1$.
However, the common definition of the reversal of a polynomial matrix  is that given in (\ref{reversal}). Hence, the reversal of a constant matrix $A_0$ is $A_0$,   therefore the partial multiplicities of   $\infty$ in a constant matrix of rank $r$ are $e_1=\dots=e_r=0$. 
In this paper we adopt the second option; i.e., we will work with the reversal of a polynomial matrix according to the definition given in (\ref{reversal}).
 \end{remark}

In the matrix pencil literature, it is also common to join together the invariant factors and the  partial multiplicities of $\infty$ in the so called homogeneous invariant  factors (see, for instance,  \cite[Ch. 2]{Fri16} or \cite[Ch. 12]{Gant66}).  This concept can  also be extended  to polynomial matrices (see \cite[Ch. 2]{Fri16},  \cite{ZaTi12}).   
In fact, if $\alpha_1(s)\mid  \dots \mid \alpha_r(s)$ are the invariant factors  of a polynomial matrix $P(s)$ and     
$e_1\leq\cdots\leq e_r$
its partial multiplicities of $\infty$,  
they can be summarized in a chain of monic homogeneous polynomials $\gamma_1(s,t) \mid \dots \mid \gamma_r(s,t),\ \gamma_i(s,t) \in \FF[s,t], \ 1\leq i\leq r$, called the  {\em homogeneous invariant  factors} of $P(s)$, of the form 
\begin{equation}\label{hif}
\gamma_i(s,t)=t^{e_i}t^{\deg(\alpha_i)}\alpha_i\left(\frac st\right), \quad 1\leq i\leq r.
\end{equation}

In turn, the homogeneous invariant factors of a polynomial matrix determine the invariant factors $\alpha_i(s)=\gamma_i(s,1), \ 1\leq i \leq r$ and the partial multiplicities of $\infty,$ $e_1 \leq\dots \leq e_r$.

We will take $\gamma_i(s,t)=1$ whenever $i<1$ and $\gamma_i(s,t)=0$  when $i>r$. As a consequence, $\alpha_i(s)=1$ and $e_i=0$ for $i<1$, and  $\alpha_i(s)=0$ for $i>r$. 
 We also agree that $e_i=+\infty$ for $i>r$.

\begin{remark}\label{remgamma1}
Notice that for a polynomial matrix $P(s)$ of degree $d$ and leading coefficient $P_d$, $\rev P(0)=P_d\neq 0$, which means that $e_1=n_1(\infty,P(s))=n_1(0,\rev P(t))=0$ (see \cite[Lemma 2.7]{DeDoVa15}),  equivalently,
\begin{equation}\label{eqgamma1}
  \gamma_1(s,0)\neq 0.
  \end{equation}
\end{remark}

In what follows we will work with the homogeneous invariant factors instead of working with the invariant factors and the infinite elementary divisors separately, for it significantly simplifies expressions. Indeed, observe that given $\alpha(s), \beta(s)\in \FF[s]$ and two non negative integers $k, \ell$,  if $\phi(s,t)=t^{k}t^{\deg(\alpha)}\alpha \left(\frac st\right)$ and $\psi(s,t)=t^{\ell}t^{\deg(\beta)}\beta\left(\frac st\right)$ then
$$\left.\begin{array}{c}
\alpha(s) \mid \beta(s), \\
k \leq \ell
\end{array}\right\} \quad \textrm{ if and only if } \quad 
\phi(s,t) \mid \psi(s,t).
$$

Unlike the case of matrix pencils, the homogeneous invariant factors and the minimal indices are not a complete set of invariants for the strict equivalence of matrix polynomials (as for matrix pencils, $P(s) \se \bar P(s)$ if  $\bar P(s)=QP(s)R, \ Q, R$ non singular) (see \cite[Subsection 3.1]{DeDoMa14}). In the
lack of better structure, we shall use this one, that is,  the eigenstructure. Given $P(s), \bar P(s)\in\efe[s]^{m\times n}$  we will write $P(s)\approx \bar  P(s)$ if they have the same  eigenstructure.

For non constant matrix pencils, the sum of the degrees of the homogeneous invariant factors  plus the sum of the minimal indices is equal to the rank. 
The following is an extension of that result to polynomial matrices given in \cite{DeDoMa14}; we state it here  in terms of the homogeneous invariant factors.

\begin{lemma}[{\rm Index Sum Theorem for Matrix Polynomials \cite[Theorem 6.5]{DeDoMa14}}]
\label{theoDeDoMa}
 Let $P(s)\in\efe[s]^{m\times n}$, $\deg(P(s))=d$, 
 $\rank(P(s))=r$.
Let $\gamma_1(s,t)\mid\cdots\mid\gamma_r(s,t)$, 
$d_1\geq\cdots\geq d_{n-r}$ and $v_1\geq\cdots\geq v_{m-r}$
be the homogeneous invariant factors, column minimal indices and row minimal indices of $P(s)$, respectively. Then, 
\begin{equation}\label{eqIST}
\sum_{i=1}^r \deg(\gamma_i)+\sum_{i=1}^{n-r} d_i+\sum_{i=1}^{m-r} v_i=rd.
\end{equation}
\end{lemma}

Notice that if $P(s)\approx \bar P(s)$ then $\rank(P(s))=\rank(\bar P(s))$ (they have the same number of invariants), and  therefore,  by Lemma \ref{theoDeDoMa}, $\deg(P(s))=\deg(\bar P(s))$.
Notice also that if $P(s), \bar P(s)\in\efe[s]^{m\times n}$ are matrix pencils,  $P(s)\approx \bar P(s)$ if and only if $P(s)\se \bar P(s)$.

Later in this paper we will frequently use Lemma \ref{theoDeDoMa}  without making any further reference to it.

\subsection{Row   completion of matrix pencils} \label{secpencilcompletion} 

In \cite{Do13} (see also \cite{DoSt19}) a solution  to Problem \ref{problem}  is given when the eigenstructure is prescribed and  the polynomial matrix is of degree one.  
We bring here the result (see Theorem \ref{theopencilcompletion} below); it  will be used later to solve the general case.  
The theorem provides a solution of the row completion problem; 
by transposition, and interchanging the row and column minimal indices, the result applies for the column completion problem.  
    To state it we need some notation and  definitions.

Given two integers $n$ and $m$, whenever $n>m$ we take  $\sum_{i=n}^{m}=0$. In the same way, if a condition is stated for $n\leq i\leq m$ with $n>m$, we understand that the condition disappears.

Let $a_1, \dots,  a_m$ be a  sequence of integers. Whenever we write  $\ba=(a_1, \dots,  a_m)$, we will understand that   $a_1\geq \dots \geq a_m$,
 and  we will take $a_i=\infty$ for $i<1$ and $a_i=-\infty$ for $i>m$. If $a_m\geq 0$, the sequence $\ba=(a_1, \dots,  a_m)$ is called a {\em partition}.  

Let   $\ba= (a_1,  \ldots, a_m)$ and $\bb= (b_1, \ldots, b_m)$ be  two sequences of integers. It is said that   $\ba$ is {\em majorized} by $\bb$ (denoted by $\ba \prec \bb$) if $\sum_{i=1}^k a_i \leq \sum_{i=1}^k b_i $ for $1 \leq k \leq m-1$ and $\sum_{i=1}^m a_i =\sum_{i=1}^m b_i$ (this is an extension to sequences of integers of the definition of majorization given for partitions in \cite{HLP88}).

We introduce next the concept of generalized majorization. 

\begin{definition}{\rm \cite[Definition 2]{DoStEJC10}}
Let $\bd = (d_1, \dots, d_m)$, $\ba=(a_1, \dots, a_s)$ and $\bg=(g_1, \dots, g_{m+s})$  be sequences of  integers.
We say that  $\bg$ is majorized by $\bd$ and $\ba$  $(\bg \prec' (\bd,\ba))$ if
\begin{equation}\label{gmaj1}
d_i\geq g_{i+s}, \quad 1\leq i\leq m,
\end{equation}
\begin{equation}\label{gmaj2}
\sum_{i=1}^{h_j}g_i-\sum_{i=1}^{h_j-j}d_i\leq \sum_{i=1}^j a_i, \quad 1\leq j\leq s,
\end{equation}
where $h_j=\min\{i\; : \; d_{i-j+1}<g_i\}$, $1\leq j\leq s$  $(d_{m+1}=-\infty)$,
\begin{equation}\label{gmaj3}
\sum_{i=1}^{m+s}g_i=\sum_{i=1}^md_i+\sum_{i=1}^sa_i.
\end{equation}
\end{definition}

\begin{remark}\label{rem_maj} \
\begin{enumerate} 

    \item \label{haches}
    The definition of $h_j$ implies that $j\leq h_j\leq m+j, 1\leq j\leq s$. 

    \item \label{gm_pc}
    In the case that $s=0$, condition (\ref{gmaj2}) disappears, and conditions (\ref{gmaj1}) and (\ref{gmaj3}) are equivalent to  $\bd=\bg$. 
On the other hand, if $m=0$ then $\bg \prec' (\bd,\ba)$ is equivalent to  $\bg \prec \ba$.

    \item \label{remg+k}
  Let $k$ be an integer, 
  $\bar \bd = (d_1+k, \dots, d_m+k)$, $\bar \ba=(a_1+k, \dots, a_s+k)$ and $\bar \bg=(g_1+k, \dots, g_{m+s}+k)$. Then
  $\bg\prec'(\bd, \ba)$ if and only if $\bar \bg\prec'(\bar \bd,\bar  \ba)$.
  \end{enumerate}
\end{remark}

The next result is in \cite{DoSt19}; we state it for non constant pencils.

\begin{theorem}\mbox{\rm \cite[Theorem 4.3]{DoSt19}}
  \label{theopencilcompletion}
Let $C(s)\in  \efe[s]^{(\bar r+p)\times (\bar r+q)}$ be a non constant matrix
pencil, $\rank(C(s))=\bar r$. Let $\bar \phi_1(s,t)\mid \dots  \mid\bar \phi_{\bar r}(s,t)$  be its homogeneous invariant factors, $\bar \bc=(\bar c_1,  \dots, \bar c_q)$  its
column minimal indices, and $\bar \bu=(\bar u_1, \dots, \bar u_p)$  its row minimal
indices, where   $\bar u_1 \geq \dots \geq \bar u_{\theta}  > \bar u_{\theta +1}= \dots = \bar u_p=0 $. Let $x$ and $y$ be non negative integers. 
Let $D(s) \in  \efe[s]^{(\bar r+p+x+y)\times (\bar r+q)}$ be a matrix pencil,  $\rank(D(s))=\bar r+x$. Let $\bar \gamma_1(s,t) \mid \dots  \mid\bar \gamma_{\bar r+x}(s,t)$
be its homogeneous invariant factors, $\bar \bd=(\bar d_1,  \dots,  \bar d_{q - x})$  its column minimal
indices, and  $\bar  \bv=(\bar v_1, \dots, \bar v_{p+y})$   its row minimal indices, where $\bar v_1 \geq  \dots  \geq  \bar v_{\bar \theta}  > \bar v_{\bar \theta +1} = \dots =\bar v_{p+y} = 0$.
There exists a pencil $A(s)$ such that
$\begin{bmatrix}C(s)\\A(s)\end{bmatrix}\se D(s)$
if and only if
\begin{equation}\label{eqinterhom}
 \bar \gamma_i(s,t)\mid \bar \phi_i(s,t)\mid \bar \gamma_{i+x+y}(s,t),\quad 1\leq i \leq \bar r,
\end{equation}
\begin{equation}\label{eqtheta}
\bar \theta \geq  \theta,
\end{equation}
\begin{equation}\label{eqcmimaj}
 \bar \bc  \prec'  (\bar \bd , \bar \ba),
\end{equation}
\begin{equation}\label{eqrmimaj}
\bar  \bv \prec'  (\bar \bu , \bar \bb),
\end{equation}
\begin{equation}\label{eqdegsum}
\sum_{i=1}^{\bar r+x}\deg(\lcm(\bar \phi_{i-x},\bar \gamma_i))\leq \sum_{i=1}^{p+y}\bar v_i-\sum_{i=1}^{p}\bar u_i+\sum_{i=1}^{\bar r+x}\deg(\bar \gamma_i),
\end{equation}
where $\bar \ba = (\bar a_1, \dots, \bar a_x )$ and $\bar \bb = (\bar b_1, \dots,  \bar b_y )$ are
\begin{equation}\label{eqdefbara}
\begin{array}{rl}
\bar a_1=&
\sum_{i=1}^{p+y}\bar v_i-\sum_{i=1}^{p}\bar u_i+\sum_{i=1}^{\bar r+x}\deg(\bar \gamma_i)-\sum_{i=1}^{\bar r+x-1}\deg(\lcm(\bar\phi_{i-x+1},\bar \gamma_i))-1,\\
\bar a_j=&
\sum_{i=1}^{\bar r+x-j+1}\deg(\lcm(\bar\phi_{i-x+j-1},\bar \gamma_i))-\sum_{i=1}^{\bar r+x-j}\deg(\lcm(\bar\phi_{i-x+j},\bar \gamma_i))-1,\\
& \hfill 2\leq j \leq x,
\end{array}
\end{equation}
\begin{equation}\label{eqdefbarb}
\begin{array}{rl}
\bar b_1=&
\sum_{i=1}^{p+y}\bar v_i-\sum_{i=1}^{p}\bar u_i+\sum_{i=1}^{\bar r+x}\deg(\bar \gamma_i)-\sum_{i=1}^{\bar r+x}\deg(\lcm(\bar\phi_{i-x-1},\bar \gamma_i)),\\
\bar b_j=&
\sum_{i=1}^{\bar r+x}\deg(\lcm(\bar\phi_{i-x-j+1},\bar \gamma_i))-\sum_{i=1}^{\bar r+x}\deg(\lcm(\bar\phi_{i-x-j},\bar \gamma_i)),\; 2\leq j \leq y.
\end{array}
\end{equation}

\end{theorem}

\begin{remark}\label{remdecreasing} \

 \begin{enumerate}
	\item \label{remdecreasing1} 
Applying  \cite[Lemmas 1 and 2]{DoStEJC10},
in Theorem \ref{theopencilcompletion}
we obtain that $\bar a_1\geq \dots\geq  \bar a_x$ and $ \bar b_1\geq \dots\geq  
\bar b_{y}\geq 0$.

\item \label{remhomified} If $C(s)$ has $\bar \alpha_1(s)\mid\cdots\mid\bar \alpha_{\bar r}(s)$ and $\bar e_1\leq\cdots\leq \bar e_{\bar r}$ as invariant factors and  partial multiplicities of $\infty$, respectively, and
$\bar \beta_1(s)\mid\cdots\mid\bar \beta_{\bar r+x}(s)$ are the invariant factors of $D(s)$ and $\bar f_1\leq\cdots\leq \bar f_{\bar r+x}$  are its partial multiplicities of $\infty$, 
then the condition (\ref{eqinterhom}) is equivalent to 
$$ \bar \beta_i(s)\mid\bar \alpha_i(s)\mid \bar \beta_{i+x+y}(s),\quad 1\leq i \leq \bar r,
$$
$$\bar  f_i\leq \bar e_i\leq \bar f_{i+x+y},\quad 1\leq i \leq \bar r.
$$
	\end{enumerate}
\end{remark}

\subsection{Frobenius companion form} \label{Frobenius}

In this subsection we state that two polynomial matrices of the same degree have the same  eigenstructure if and only if their first Frobenius companion forms are strictly equivalent.

\begin{definition}\label{defFrpbform}
Let $P(s)=P_ds^d+P_{d-1}s^{d-1}+\cdots+P_1s+P_0\in\efe[s]^{m\times n}$, $P_d\neq 0$, $d\geq 1$. The {\em first Frobenius companion form} of $P(s)$ is the $(m+(d-1)n)\times dn$ pencil $C_P(s)=sX_1+Y_1$ with
$$
X_1=\begin{bmatrix}P_d&&&\\&I_n&&\\&&\ddots&\\&&&I_n\end{bmatrix}\,\text{ and }\,
Y_1=\begin{bmatrix}P_{d-1}&P_{d-2}&\cdots&P_0\\-I_n&0&\cdots&0\\&\ddots&\ddots&\vdots\\0&&-I_n&0\end{bmatrix}.
$$
\end{definition}	

Notice  that when $d=1$, $C_P(s)=P(s)$.
The following lemma is a consequence of \cite[Theorem 5.3 and Theorem 4.1]{DeDoMa14}.

\begin{lemma}\label{lem_lin}
	Let $P(s)\in\efe[s]^{m\times n}$  be a polynomial matrix of degree $d\geq 1$ and let $C_P(s)$ be its first Frobenius companion form.  Then,
	\begin{enumerate}
	\item 
          If $\phi_1(s,t), \dots, \phi_{r}(s,t)$  are the  homogeneous invariant factors of $P(s)$ then
          the  homogeneous invariant factors of the pencil $C_P(s)$ are
 $1, \stackrel{(d-1)n}{\dots}, 1, \phi_1(s,t), \dots,  \phi_{r}(s,t)$.

	\item 
	If $c_1\geq\cdots\geq c_{n-r}$ are the column minimal indices of $P(s)$ then $c_1+d-1\geq\cdots\geq c_{n-r}+d-1$ are the column minimal indices of $C_P(s)$.
\item 
	If $u_1\geq\cdots\geq u_{m-r}$ are the row minimal indices of $P(s)$ then $u_1\geq\cdots\geq u_{m-r}$ are the row minimal indices of $C_P(s)$.
        \end{enumerate}       
\end{lemma}

As an immediate consequence of Lemma \ref{lem_lin} we obtain  the next corollary.

\begin{corollary}\label{corcompofse}
Let $P(s), \bar P(s)\in\efe[s]^{m\times n}$ be polynomial matrices satisfying that $\deg(P(s))=\deg(\bar P(s))\ge 1$, and let $C_P(s),C_{\bar P}(s) $ be their first Frobenius companion forms, respectively.  Then,
$P(s)\approx \bar P(s)$ if and only if $C_P(s)\se C_{\bar P}(s)$.	
  \end{corollary}

\section{Polynomial matrices with prescribed eigenstructure} \label{secpolynomialprescreig} 

The aim of this section is to show a characterization of    the existence of a polynomial matrix of degree $d$  with prescribed eigenstructure over arbitrary fields.
For $d\geq 1$, the result was obtained in 
\cite[Theorem 3.3]{DeDoVa15} 
but in that paper the proof  of the sufficiency  is valid only 
when $\efe$ is an 
infinite field. In Theorem \ref{theoexistenceDeDoVa152} we remove this restriction. We also include the case $d=0$.

  For matrix pencils the sufficiency follows from the Kronecker canonical form.
For $d>1$ we use the first Frobenius companion form for transforming the problem into a row completion problem of the subpencil formed by the last $(d-1)n$ rows of this companion form (see the pencil $C(s)$ in (3.1)). 
 With this idea in mind, by Theorem  \ref{theopencilcompletion}, the conditions of \cite[Theorem 3.3]{DeDoVa15} allow us to complete $C(s)$ up to a pencil  with the desired eigenstructure.

\begin{theorem}\label{theoexistenceDeDoVa152} 
Let $m$, $n$, $r\leq \min\{m,n \}$ be  positive integers and $d$ a non negative integer.
Let
$\gamma_1(s,t)\mid \dots \mid \gamma_r(s,t)$ be monic homogeneous polynomials  with coefficients in an arbitrary field   $\efe$.
Let $(d_1, \ldots, d_{n-r})$, $(v_1, \ldots, v_{m-r})$
be partitions. 
Then there exists $P(s)\in \efe[s]^{m\times n}$, $\rank(P(s))= r$,
$\deg(P(s))= d$,  with homogeneous invariant factors $\gamma_1(s,t), \dots, \gamma_r(s,t)$ and column and row minimal indices
$d_1, \dots, d_{n-r}$ and $v_1, \dots, v_{m-r}$, respectively, if and only if
(\ref{eqgamma1}) and (\ref{eqIST}) hold.
\end{theorem}

\noindent
{\bf Proof}.
The proof of the necessity in \cite[Theorem 3.3]{DeDoVa15} is valid for arbitrary fields. Observe that it also holds in the case that $d=0$.
Therefore, we only need to  prove the  sufficiency.
Assume  that  (\ref{eqgamma1}) and (\ref{eqIST})
are satisfied.

 For $d=0$,
  from (\ref{eqIST}) we obtain $\gamma_1(s,t)=\dots =\gamma_r(s,t)=1$ and 
$d_1=\dots = d_{n-r}=v_1= \dots =v_{m-r}=0$.  Then, any matrix $P_0\in \FF^{m\times n}$ of rank $r$ has the desired invariants.

For $d=1$, using  the Kronecker canonical form of matrix pencils,   we can build a matrix polynomial $P(s)\in \efe[s]^{m\times n}$ of degree 1 and rank $r$
having   $\gamma_1(s,t), \dots, \gamma_r(s,t)$ as homogeneous invariant factors and  $d_1, \dots, d_{n-r}$ and $v_1, \dots, v_{m-r}$ as column and row minimal indices, respectively.
Indeed, let
$\alpha_i(s)=\gamma_i(s,1)$, $1\leq i \leq r$, and let
$e_1\leq \dots\leq  e_r$ be defined by (\ref{hif}). Assume that 
$\alpha_1(s)=\dots=\alpha_w(s)=1\neq \alpha_{w+1}(s)$ ($0\leq w\leq r$), $e_1= \dots =e_q=0 <e_{q+1}$ ($1\leq q\leq r$), $d_1\geq \dots \geq d_\rho>0$ ($0\leq \rho\leq n-r$),  $v_1\geq \dots \geq v_\theta>0$ ($0\leq \theta\leq m-r$). Notice that  because of (\ref{eqgamma1}), $q\geq 1$.
Let
$$
P(s)=\left[\begin{matrix}\diag(C(s), N(s),L(s), R(s) )&0\\0&0\end{matrix}\right]\in \FF[s]^{m\times n},
$$
where $$C(s)=\diag(C_{\alpha_{w+1}}(s), \dots, C_{\alpha_r}(s)),\quad 
N(s)=\diag(N_{e_{q+1}}(s), \ldots, N_{e_r}(s)),$$ 
$$L(s)=\diag(L_{d_1}(s), \ldots, L_{d_\rho}(s)),  \mbox{ and }R(s)=\diag(R_{v_1}(s), \ldots, R_{v_\theta}(s)).$$  
Here the block matrix $C_{\alpha}(s) \in \FF[s]^{\deg(\alpha)\times \deg(\alpha)}$ is the first Frobenius companion form of  $\alpha(s)$,
$$
N_k(s)={\scriptsize 
\begin{bmatrix}
1 & s   &  & \\
  &  \ddots & \ddots &    \\
 &  &  \ddots & s \\
  &  &  & 1
\end{bmatrix}
}
\in \FF[s]^{k\times k},\quad 
L_k(s)= {\scriptsize
\begin{bmatrix}
s & 1 &   &  \\
 &    \ddots & \ddots &    \\
 & &  s & 1
 \end{bmatrix}
 }
\in \FF[s]^{k \times (k+1)},
$$
and  $R_k(s)=L_k(s)^T\in \FF[s]^{(k+1)\times k}.$
If $\deg(\alpha_r(s))>0$ or $d_1>0$ or $v_1>0$ then $\deg(P(s))=1$. Otherwise, by (\ref{eqgamma1}) and (\ref{eqIST}), $\sum_{i=2}^r e_i=r$. Thus, $e_r>1$ and $\deg(P(s))=1$. Now, by (\ref{eqIST}), $\rank(P(s))= r$, and $P(s)$ is the desired pencil.

Let $d>1$ and
\begin{equation}\label{eqlastrows}
C(s)=\begin{bmatrix}-I_n&sI_n&0&\cdots&0&0\\
0&-I_n&sI_n&\cdots&0&0\\
\vdots &\vdots &\ddots&\ddots&\vdots&\vdots \\0&0&0&\dots&sI_n&0\\0&0&0&\dots&-I_n&sI_n
\end{bmatrix}\in \efe[s]^{(d-1)n\times dn}.
\end{equation}

Then, $\bar  r=\rank(C(s))=(d-1)n$, the
homogeneous invariant factors of $C(s)$ are
$\bar \phi_1(s,t) =\dots =\bar \phi_{\bar r}(s,t)=1$, 
and the columns of the matrix $\begin{bmatrix}s^{d-1}I_n\\s^{d-2}I_n\\\vdots\\sI_n\\I_n
\end{bmatrix}$ form 
a right minimal basis for $C(s)$; hence,  denoting the column minimal indices of $C(s)$  by $\bar \bc=(\bar c_1, \dots, \bar c_n)$, we have
$\bar c_i=d-1$, $1\leq i \leq n$.
As $C(s)$ does  not have row minimal indices, we put $\bar \bu=\emptyset$. We also take $\theta=0$.

We introduce now a collection of homogeneous polynomials and two partitions of integers intended to be the Kronecker invariants of an $(m+(d-1)n) \times dn$ pencil. 
Define
$$\begin{array}{rl}
          \bar \gamma_i(s,t)=&1, \quad 1\leq i \leq \bar r, \\
          \bar \gamma_{\bar r+i}(s,t)=&\gamma_i (s, t),  \quad 1\leq i \leq r, \\
\end{array}$$
$$
\bar \bd=(\bar d_1, \dots, \bar d_{n-r}), \quad 
\bar \bv=(\bar v_1, \dots, \bar v_{m-r}),$$
where
$$\bar d_i=d_i+d-1, \quad  1\leq i \leq n-r, $$  $$\bar v_i=v_i, \quad  1\leq i \leq m-r.$$
Observe that $n-r=dn-(\bar r+r)$ and $m-r=m+(d-1)n-(\bar r+r)$,
and taking into account (\ref{eqIST}),
$$\sum_{i=1}^{\bar r+r} \deg(\bar \gamma_i)+\sum_{i=1}^{n-r} \bar d_i+\sum_{i=1}^{m-r} \bar v_i=rd+(n-r)(d-1)=\bar r+r.$$
  Applying the result for $d=1$,
there is a pencil $D(s)\in \efe[s]^{(m+(d-1)n)\times dn}$, $\deg(D(s))=1$, $\rank(D(s))=\bar r+r$, with
homogeneous invariant factors $\bar \gamma_1(s,t), \dots, \bar \gamma_{\bar r+r}(s,t) $, column minimal indices
$\bar d_1, \dots, \bar d_{n-r}$, and row minimal indices $\bar v_1, \dots, \bar v_{m-r}$.

Let us see that the invariants of the pencils $C(s)$ and $D(s)$ satisfy the conditions (\ref{eqinterhom})-(\ref{eqdegsum}) of Theorem  \ref{theopencilcompletion}.
Take $x=r$, $y=m-r$, and $\bar \theta=\#\{i\; : \; \bar v_i>0\}$. Then 
(\ref{eqinterhom}) and (\ref{eqtheta}) hold.
 Since
$
\sum_{i=1}^{\bar r+r}\deg(\lcm(\bar \phi_{i-r},\bar \gamma_i))=\sum_{i=1}^{\bar r+r}\deg(\bar \gamma_i)$ and
 $\sum_{i=1}^{m-r}\bar v_i\geq 0$, condition (\ref{eqdegsum}) holds.

Let $\bar \ba = (\bar a_1, \dots, \bar a_r )$ and $\bar \bb = (\bar b_1, \dots,  \bar b_{m-r} )$
be defined  as in (\ref{eqdefbara}) and (\ref{eqdefbarb}), respectively. 
Trivially  $\bar \bv\prec\bar \bb$, which is equivalent to (\ref{eqrmimaj}) (see Remark \ref{rem_maj}.\ref{gm_pc}). 
Finally,
$$
\bar d_{i}=d_{i}+d-1\geq d-1=\bar c_{i+r}, \quad 1\leq i \leq n-r,
$$
$$
\begin{array}{rl}
\sum_{i=1}^{n-r}\bar d_i+\sum_{i=1}^{r}\bar a_i=&\sum_{i=1}^{n-r}d_i+(n-r)(d-1)+\sum_{i=1}^{m-r} v_i+\sum_{i=1}^{r}\deg(\gamma_i)-r\\=&
rd+(n-r)(d-1)-r=n(d-1)=\sum_{i=1}^{n}\bar c_i.\end{array}
$$
Let $j\in \{1, \dots, r\}$. Let $h_j=\min\{i\; : \; \bar d_{i-j+1}<\bar c_i\}$. By Remark \ref{rem_maj}.\ref{haches}, $j\leq h_j\leq n-r+j$.
For $j\leq i \leq (n-r)+j-1$,
we have
$$
\bar d_{i-j+1}=d_{i-j+1}+d-1\geq d-1=\bar c_i;
$$
hence, $h_j=n-r+j$ and
$$\begin{array}{l}\sum_{i=1}^{j}\bar a_i+\sum_{i=1}^{h_j-j}\bar d_i-\sum_{i=1}^{h_j}\bar c_i
=\sum_{i=1}^{m-r}v_i+\sum_{i=r-j+1}^{r}\deg(\gamma_i)-j\\
+\sum_{i=1}^{n-r}d_i+(n-r)(d-1)-(n-r+j)(d-1)\\
=\sum_{i=1}^{r}\deg(\gamma_i)+\sum_{i=1}^{n-r}d_i+\sum_{i=1}^{m-r}v_i-\sum_{i=1}^{r-j}\deg(\gamma_i)-jd.\end{array}
$$
From (\ref{eqIST}), we obtain
$
\sum_{i=1}^{j}\bar a_i+\sum_{i=1}^{h_j-j}\bar d_i-\sum_{i=1}^{h_j}\bar c_i=(r-j)d-\sum_{i=1}^{r-j}\deg(\gamma_i)$ and  $\sum_{i=1}^{r}\deg(\gamma_i)\leq rd$. As 
$\deg(\gamma_1)\leq \dots \leq \deg(\gamma_r)$, it follows that 
 $\sum_{i=1}^{r-j}\deg(\gamma_i)\leq (r-j)d$. 
Therefore, $(r-j)d-\sum_{i=1}^{r-j}\deg(\gamma_i)\geq 0$, and 
condition (\ref{eqcmimaj}) holds.
By Theorem \ref{theopencilcompletion}, there exists a pencil $W(s)\in \efe[s]^{m\times dn}$ such that
$\begin{bmatrix}W(s)\\C(s)\end{bmatrix}$ has $\bar \gamma_1(s,t), \dots, \bar \gamma_{\bar r+r}(s,t) $ as  homogeneous invariant factors, the integers $\bar d_1, \dots, \bar d_{n-r}$ as column minimal indices, and $\bar v_1, \dots, \bar v_{m-r}$ as row minimal indices.
Let
$$
\begin{bmatrix}W(s)\\C(s)\end{bmatrix}=\begin{bmatrix}
W_{d-1}(s)&W_{d-2}(s)&W_{d-3}(s)&\cdots&W_{1}(s)&W_{0}(s)\\
-I_n&sI_n&0&\cdots&0&0\\
0&-I_n&sI_n&\cdots&0&0\\
\vdots &\vdots &\ddots&\ddots&\vdots&\vdots\\0&0&0&\dots&sI_n&0 \\0&0&0&\dots&-I_n&sI_n\end{bmatrix},
$$
and let $W_i(s)=sW_{i,1}+W_{i,0}\in \efe[s]^{m\times n}$, $0\leq i \leq d-1$. There exists  $U\in \Gl_{m+(d-1)n}(\FF)$ such that
$$
\begin{bmatrix}W(s)\\C(s)\end{bmatrix}\se U\begin{bmatrix}W(s)\\C(s)\end{bmatrix}=
\begin{bmatrix}
sP_d+P_{d-1}&P_{d-2}&P_{d-3}&\cdots&P_{1}&P_{0}\\
-I_n&sI_n&0&\cdots&0&0\\
0&-I_n&sI_n&\cdots&0&0\\
\vdots &\vdots &\ddots&\ddots&\vdots&\vdots\\0&0&0&\dots&sI_n&0 \\0&0&0&\dots&-I_n&sI_n\end{bmatrix},
$$
where
$P_d=W_{d-1,1}$, $P_0=W_{0,0}$,  and $P_i=W_{i,0}+W_{i-1,1}$, $1\leq i \leq d-1$.

Let us denote $K(s)=U\begin{bmatrix}W(s)\\C(s)\end{bmatrix}$. From (\ref{eqgamma1}) we derive  $\bar \gamma_{\bar r+1}(s,0)\neq 0$. Thus, $n_i(\infty,K(s))=n_i(0,\rev K(t))=0$, $1\leq i \leq \bar r+1$. Observe that
$$
\rev K(t)=\begin{bmatrix} P_d&0\\0&I_{\bar r}\end{bmatrix}+t
\begin{bmatrix} \begin{matrix}P_{d-1}&\cdots&P_1\end{matrix}&P_0\\-I_{\bar r}&0\end{bmatrix},
$$
and $\rev K(0)=\begin{bmatrix} P_d&0\\0&I_{\bar r}\end{bmatrix}$. As $n_{\bar r +1}(0,\rev K(t))=0$, the $(\bar r +1)$-th invariant factor of $\rev K(t)$ is not a multiple of $t$, therefore, $\rank(\rev K(0))\geq \bar r+1$ and $P_d\neq 0$.

Let $P(s)=P_ds^d+P_{d-1}s^{d-1}+\cdots+ P_1s+P_0 \in \efe[s]^{m\times n}$. 
By Lemma \ref{lem_lin},
$P(s)$ has  homogeneous invariant factors $\gamma_1(s,t), \dots, \gamma_{r}(s,t)$, 
column minimal indices
$d_1, \dots, d_{n-r}$, and row minimal indices $v_1, \dots, v_{m-r}$.
\hfill $\Box$

\section{Row (column) completion  of polynomial matrices of given degree}
\label{secpolcompl}

The next proposition  turns a row completion problem of polynomial matrices into a row completion problem of matrix pencils.  As a solution to the latter is known, out of the solution of the problem for matrix pencils we are able to find a solution to the problem for polynomial matrices.

\begin{proposition}\label{propeqprob}
  Let $P(s)\in\efe[s]^{m\times n}$ and $Q(s)\in\efe[s]^{(m+z)\times n}$ be polynomial matrices such that  $\deg(P(s))=\deg (Q(s))=d\ge 1$,
  and let $C_P(s),C_{Q}(s) $ be their first Frobenius companion forms, respectively.
  Then, there exists $W(s)\in\efe[s]^{z\times n}$ such that $\deg(W(s))\leq d$ and  $\begin{bmatrix}P(s)\\W(s)\end{bmatrix}\approx Q(s)$ if and only if there exists a matrix pencil
   $A(s)\in\efe[s]^{z\times dn}$ such that $\begin{bmatrix}C_P(s)\\A(s)\end{bmatrix}\se C_Q(s)$.
  \end{proposition}

\noindent
{\bf Proof}.
Let $P(s)=P_ds^d+P_{d-1}s^{d-1}+\cdots+P_1s+P_0$, $P_d\neq 0$.
Assume that there exists $W(s)=W_ds^d+W_{d-1}s^{d-1}+\cdots+W_1s+W_0\in\efe[s]^{z\times n}$, $\deg(W(s))\leq d$, such that $\bar P(s)=\begin{bmatrix}P(s)\\W(s)\end{bmatrix}\approx Q(s)$.
Observe that $\deg(\bar P(s))=d$, and  let 
$$
C_{\bar P}(s)=\left[
\begin{smallmatrix}sP_d+P_{d-1}&P_{d-2}&\dots &P_1&P_0\\sW_d+W_{d-1}&W_{d-2}&\dots &W_1&W_0\\
  -I_n&sI_n&\ddots &0&0\\
  0&-I_n&\ddots &0&0\\
  \vdots &\vdots &\ddots &\vdots &\vdots\\
  0&0&\dots &-I_n&sI_n
\end{smallmatrix}\right]\in \efe[s]^{(m+z+(d-1)n)\times dn},
$$
be the first Frobenius companion form of $\bar P(s)$. By Corollary \ref{corcompofse},
$C_{\bar P}(s)\se C_Q(s)$.
Let 
$A(s)=\begin{bmatrix}sW_d+W_{d-1}&W_{d-2}&\dots &W_1&W_0\end{bmatrix}\in \efe[s]^{z\times dn}$. Then,
$C_Q(s)\se C_{\bar P}(s)\se \begin{bmatrix}C_P(s)\\A(s)\end{bmatrix}.$

Conversely, 
assume that there exists  a matrix pencil
$$A(s)=\begin{bmatrix}sA_{d-1,1}+A_{d-1,0}&sA_{d-2,1}+A_{d-2,0}&
&\dots &sA_{0,1}+A_{0,0}
\end{bmatrix}\in \efe[s]^{z\times dn},$$ such that $\begin{bmatrix}C_P(s)\\A(s)\end{bmatrix}\se C_Q(s)\in \efe[s]^{(m+(d-1)n+z)\times dn}$. Let 
$$
\hat C(s)=\begin{bmatrix}C_P(s)\\A(s)\end{bmatrix}=\left[
\begin{smallmatrix}sP_d+P_{d-1}&P_{d-2}&\dots &P_1&P_0\\
  -I_n&sI_n&\ddots &0&0\\
  0&-I_n&\ddots &0&0\\
  \vdots &\vdots &\ddots &\vdots &\vdots\\
  0&0&\dots &-I_n&sI_n\\
  sA_{d-1,1}+A_{d-1,0}&sA_{d-2,1}+A_{d-2,0}
&\dots &sA_{1,1}+A_{1,0}&sA_{0,1}+A_{0,0}
\end{smallmatrix}\right].
$$
There exists $U\in \Gl_{m+z+(d-1)n}(\efe)$  such that 
$$U\hat C(s)=\left[
\begin{smallmatrix}sP_d+P_{d-1}&P_{d-2}&\dots &P_1&P_0\\
  sA_{d-1,1}+A'_{d-1}&A'_{d-2}&\dots &A'_1&A'_0\\
  -I_n&sI_n&\ddots &0&0\\
  0&-I_n&\ddots &0&0\\
  \vdots &\vdots &\ddots &\vdots &\vdots\\
  0&0&\dots &-I_n&sI_n
\end{smallmatrix}\right]\in \efe[s]^{(m+z+(d-1)n)\times dn},
$$
where $A'_i=A_{i,0}+A_{i-1,1}$, $0\leq i \leq d-1$ ($A_{-1,1}=0$).
Let $W(s)=s^dA_{d-1,1}+s^{d-1}A'_{d-1}+\dots+sA'_1+A'_0\in \efe[s]^{z\times n}$ and
$\bar P(s)=\begin{bmatrix}P(s)\\W(s)\end{bmatrix}$. 
 Then
$C_{\bar P}(s)=U\hat C(s)\se C_Q(s)$.  By Corollary \ref{corcompofse},
$\bar P(s)\approx Q(s)$.
\hfill $\Box$

\subsection{Prescription of the whole  eigenstructure}\label{subsec_eigen}

The next theorem  generalizes  Theorem \ref{theopencilcompletion} to polynomial matrices. It contains a solution to Problem \ref{problem} when the whole eigenstructure is prescribed. In subsequent subsections different particular cases of the problem  are solved when only some of the  invariants are prescribed. Our target is to analyze all of the possible cases, but as  there are many possibilities, we will present here only some of them. In a future paper we will accomplish the study of the remaining cases. 
\begin{theorem}
	\label{theoprescr4}
	Let $P(s)\in\efe[s]^{m\times n}$ be a polynomial matrix, $\deg(P(s))=d$, $\rank (P(s))=r$.
	Let $\phi_1(s,t)\mid \cdots\mid\phi_r(s,t)$ be its
	homogeneous  invariant factors, 
	$\bc=(c_1,  \dots,  c_{n-r})$  its
	column minimal indices, and $\bu=(u_1, \dots,  u_{m-r})$  its row minimal
	indices, where   $u_1 \geq \dots \geq u_{\eta}  >  u_{\eta +1}= \dots = u_{m-r}=0 $.
	
	Let $z$ and $x$ be integers such that $0\leq x\leq \min\{z, n-r\}$.  
	Let $\gamma_1(s,t)\mid\cdots\mid\gamma_{r+x}(s,t)$ be monic homogeneous polynomials, and
	$\bd=(d_1, \dots, d_{n-r-x})$
	 and 
	$\bv=(v_1, \dots, v_{m+z-r-x})$ two partitions, where 
	$v_1 \geq  \dots  \geq  v_{\bar \eta}  > v_{\bar \eta +1} = \dots =v_{m+z-r-x} = 0$.
	There exists a polynomial matrix  $W(s)\in \efe[s]^{z\times n}$ such that $\deg(W(s))\leq d$, $\rank \left(\begin{bmatrix}P(s)\\W(s)\end{bmatrix}\right) =r+x$, and 
	$\begin{bmatrix}P(s)\\W(s)\end{bmatrix}$ has $\gamma_1(s,t)\mid\cdots\mid\gamma_{r+x}(s,t)$ as homogeneous invariant factors,
	$d_1, \dots, d_{n-r-x}$ as column minimal indices
	and 
	$v_1, \dots, v_{m+z-r-x}$ as row minimal indices  
	if and only if 
	\begin{equation}\label{eqinterfipolhom}
		\gamma_i(s,t)\mid \phi_i(s,t)\mid \gamma_{i+z}(s,t),\quad 1\leq i \leq r,
	\end{equation}
	\begin{equation}\label{eqetapol}\bar \eta \geq  \eta,\end{equation}
	\begin{equation}\label{eqcmimajpol}  \bc \prec'  (\bd , \ba),\end{equation}
	\begin{equation}\label{eqrmimajpol}\bv \prec'  (\bu , \bb),\end{equation}

\begin{equation}\label{eqdegsumpol}
 \begin{aligned}
		\sum_{i=1}^{r+x}\deg(\lcm(\phi_{i-x},\gamma_i))
		\leq \sum_{i=1}^{m+z-r-x}v_i-\sum_{i=1}^{m-r}u_i+\sum_{i=1}^{r+x}\deg(\gamma_i),
  \\
\mbox{ with equality when $x=0$,}
\end{aligned}
\end{equation}
	where $\ba = (a_1, \dots, a_x )$ and $\bb = (b_1, \dots,  b_{z-x} )$ are
	\begin{equation}\label{eqdefa}
		\begin{array}{rl}
			a_1=&
			\sum_{i=1}^{m+z-r-x}v_i-\sum_{i=1}^{m-r} u_i\\&+\sum_{i=1}^{r+x}\deg( \gamma_i)-
			\sum_{i=1}^{  r+x-1}\deg(\lcm(  \phi_{i-x+1},  \gamma_i))-d,\\
			a_j=&
			\sum_{i=1}^{  r+x-j+1}\deg(\lcm(  \phi_{i-x+j-1},  \gamma_i))
			-
			\sum_{i=1}^{  r+x-j}\deg(\lcm(  \phi_{i-x+j},  \gamma_i))
			-d,\\ & \hfill 2\leq j \leq x,
		\end{array}
	\end{equation}
	\begin{equation}\label{eqdefb}
		\begin{array}{rl}
			b_1=&
			\sum_{i=1}^{m+z-r-x}  v_i-\sum_{i=1}^{m-r}  u_i+\sum_{i=1}^{  r+x}\deg(  \gamma_i)-
			\sum_{i=1}^{  r+x}\deg(\lcm(  \phi_{i-x-1},  \gamma_i)),
			\\
			b_j=&
			\sum_{i=1}^{  r+x}\deg(\lcm(  \phi_{i-x-j+1},  \gamma_i)-
			\sum_{i=1}^{  r+x}\deg(\lcm(  \phi_{i-x-j},  \gamma_i)),\\& \hfill 2\leq j \leq z-x.
		\end{array}
	\end{equation}

\end{theorem}

\noindent
{\bf Proof}.

 First notice that $\ba$ and $\bb$ are well defined (see Remark \ref{remdecreasing}.\ref{remdecreasing1}).

If $d=0$, 
then $\phi_1(s,t)=\cdots=\phi_r(s,t)=1$, 
$c_1= \dots=c_{n-r}=0$, $u_{1}= \dots = u_{m-r}=0$.

If  $d\geq 1$, take 
$\bar r=(d-1)n+r$, $ y=z-x$, $p=m-r=m+(d-1)n-\bar r$,  $q=n-r=dn-\bar r$,
and let $C_P(s)$ be the first Frobenius companion form of $P(s)$. Then
$C_P(s)\in \efe[s]^{(\bar r+p)\times (\bar r+q)}$ and
$\rank (C_P(s))= \bar r$. 
Let $\bar \phi_1(s,t)\mid\cdots\mid\bar \phi_{\bar r}(s,t)$, 
$\bar c_1 \geq  \dots   \geq \bar c_{q} $ 
and $\bar u_1 \geq  \dots  \geq  \bar u_{ \theta}  > \bar u_{\theta +1} = \dots =\bar u_{p} = 0$
be the homogeneous invariant factors, column minimal indices and row minimal indices of $C_P(s)$, respectively,  and let
$\bar \bc=(\bar c_1, \dots, \bar c_{q})$, and $\bar \bu=(\bar u_1,  \dots, \bar u_{p})$.
By Lemma \ref{lem_lin},
$$
  \begin{array}{l}
\bar\phi_i(s,t)=1,  \quad 1\leq i \leq (d-1)n,\\
\bar\phi_{i+(d-1)n}(s,t)=\phi_{i}(s,t), \quad   1\leq i\leq r,\\
\bar c_i=c_i+d-1, \quad 1\leq i \leq n-r=q,\\
\theta= \eta, \\
\bar u_i=u_i, \quad 1\leq i \leq m-r=p.
  \end{array}
  $$
Assume that
  there exists   $W(s)\in \efe[s]^{z\times n}$ such that $\deg(W(s))\leq d$, $\rank \left(\begin{bmatrix}P(s)\\W(s)\end{bmatrix}\right) =r+x$, and 
	$\begin{bmatrix}P(s)\\W(s)\end{bmatrix}$ has $\gamma_1(s,t)\mid\cdots\mid\gamma_{r+x}(s,t)$ as homogeneous invariant factors,
	$d_1, \dots, d_{n-r-x}$ as column minimal indices
	and 
	$v_1, \dots, v_{m+z-r-x}$ as row minimal indices. Put        
$Q(s)=\begin{bmatrix}P(s)\\W(s)\end{bmatrix}$, $\bd=(d_1, \dots,  d_{q-x})$ and  $\bv=( v_1,  \dots, v_{p+y})$.

If $d=0$, then 	$v_1, \dots, v_{m+z-r-x}=0$; hence
$a_1=\dots=a_x=0$, $b_1=\dots=b_{z-x}=0$ and 
(\ref{eqinterfipolhom})-(\ref{eqdegsumpol})  trivially hold.

If $d\geq 1$,  
let  $C_Q(s)\in \efe[s]^{(\bar r+p+x+y)\times (\bar r+q)}$  be the first Frobenius companion form of $Q(s)$. Then
 $\rank (C_Q(s))= \bar r+x$.
By Proposition \ref{propeqprob}   
there exists a  matrix  pencil $A(s)\in \efe[s]^{z\times n}$ such that
$\begin{bmatrix}C_P(s)\\A(s)\end{bmatrix}\se C_Q(s)$.

Let 
$\bar \gamma_1(s,t)\mid\cdots\mid\bar \gamma_{\bar r+x}(s,t)$, 
$\bar d_1 \geq  \dots   \geq \bar d_{q-x}$ 
and $\bar v_1 \geq  \dots  \geq  \bar v_{\bar \theta}  > \bar v_{\bar \theta +1} = \dots =\bar v_{p+y} = 0$
be the homogeneous invariant factors, column minimal indices and row minimal indices of $C_Q(s)$, respectively, and 
let  
$\bar \bd=(\bar d_1, \dots, \bar d_{q-x})$ and  $\bar \bv=(\bar v_1,  \dots, \bar v_{p+y})$.
By Lemma \ref{lem_lin},
$$
  \begin{array}{l}
\bar\gamma_i(s,t)=1,  \quad 1\leq i \leq (d-1)n,\\
\bar\gamma_{i+(d-1)n}(s,t)=\gamma_{i}(s,t),  \quad 1\leq i\leq r+x,\\
\bar d_i=d_i+d-1, \quad 1\leq i \leq n- r -x=q-x,\\
\bar \theta=\bar \eta,\\
\bar v_i=v_i, \quad 1\leq i \leq m+z-r-x=p+y.\end{array}
$$
By Theorem \ref{theopencilcompletion},  (\ref{eqinterhom})-(\ref{eqdegsum}) hold,
where $\bar \ba = (\bar a_1, \dots, \bar a_x )$ and $\bar \bb = (\bar b_1, \dots,  \bar b_y )$ are  defined as in (\ref{eqdefbara})  and (\ref{eqdefbarb}), respectively.
It is easy to see that
$\bar a_j=a_j+(d-1)$ for $1\leq j \leq x$ and $\bar b_j=b_j$ for $1\leq j \leq y$.
Then (\ref{eqinterhom})-(\ref{eqdegsum})  are equivalent to (\ref{eqinterfipolhom})-(\ref{eqdegsumpol}) (see Remark \ref{rem_maj}.\ref{remg+k}  for the equivalence between (\ref{eqcmimaj})  and
(\ref{eqcmimajpol})).

Assume now that 
(\ref{eqinterfipolhom})-(\ref{eqdegsumpol}) hold.
Notice that  from (\ref{eqdefa}), or from  (\ref{eqdegsumpol}) for $x=0$, we obtain
$$\sum_{i=1}^x a_i=	\sum _{i=1}^{m+z-r-x}v_i-\sum_{i=1}^{m-r} u_i+\sum_{i=1}^{r+x}\deg( \gamma_i)-
\sum_{i=1}^{  r}\deg(\lcm(  \phi_{i},  \gamma_i))-xd.
$$

From (\ref{eqinterfipolhom}) we have $\sum_{i=1}^{  r}\deg(\lcm(  \phi_{i},  \gamma_i))=\sum_{i=1}^{  r}\deg(\phi_{i})$ and, 
taking into account (\ref{eqcmimajpol}) we get
$$
\sum_{i=1}^{r+x}\deg( \gamma_i)+\sum_{i=1}^{n-r-x} d_i+\sum_{i=1}^{m+z-r-x}v_i=\sum_{i=1}^{n-r} c_i+\sum_{i=1}^{m-r} u_i+\sum_{i=1}^{  r}\deg(\phi_{i})+xd=(r+x)d.
$$

Moreover,  because of $\phi_1(s,0)\neq 0$ and (\ref{eqinterfipolhom}), we obtain $\gamma_1(s,0)\neq 0$. By Theorem \ref{theoexistenceDeDoVa152}, there exists $Q(s)\in \efe[s]^{(m+z)\times n}$, $\rank(Q(s))= r+x$,
$\deg(Q(s))= d$,  with homogeneous invariant factors $\gamma_1(s,t), \dots, \gamma_{r+x}(s,t)$ and column and row minimal indices
$d_1, \dots, d_{n-r-x}$ and $v_1, \dots, v_{m+z-r-x}$, respectively. 

If $d=0$, the $\gamma_1(s,t)=\cdots=\gamma_{r+x}(s,t)=1$,
	$d_1=\dots= d_{n-r-x}=0$ 
	and 
	$v_1, \dots, v_{m+z-r-x}=0$; hence
$a_1=\dots=a_x=0$, $b_1=\dots=b_{z-x}=0$.
Choosing $W\in \efe^{z\times n}$  such that 
	$\rank \begin{bmatrix}P\\W\end{bmatrix}=r+x$, the matrix $\begin{bmatrix}P\\W\end{bmatrix}$ has the prescribed invariants.

    If $d\geq 1$,
let $C_Q(s)\in \efe[s]^{(\bar r+p+x+y)\times (\bar r+q)}$ be the first Frobenius companion form of $Q(s)$. Then
  $\rank (C_Q(s))= \bar r+x$.
Let 
$\bar \gamma_1(s,t)\mid\cdots\mid\bar \gamma_{\bar r+x}(s,t)$, 
$\bar d_1 \geq  \dots   \geq \bar d_{q-x}$ 
and $\bar v_1 \geq  \dots  \geq  \bar v_{\bar \theta}  > \bar v_{\bar \theta +1} = \dots =\bar v_{p+y} = 0$
be the homogeneous invariant factors, column minimal indices and row minimal indices of $C_Q(s)$, respectively, and 
let  
$\bar \bd=(\bar d_1, \dots, \bar d_{q-x})$ and  $\bar \bv=(\bar v_1,  \dots, \bar v_{p+y})$.
As in the proof of the necessity, (\ref{eqinterfipolhom})-(\ref{eqdegsumpol}) are equivalent to
(\ref{eqinterhom})-(\ref{eqdegsum}).
By Theorem \ref{theopencilcompletion}, there exists a matrix pencil
$A(s) \in \FF[s]^{z\times dn}$  such that 
that $\begin{bmatrix}C_P(s)\\A(s)\end{bmatrix} \se Q(s)$.
By Proposition \ref{propeqprob}   
there exists a polynomial matrix  $W(s)\in \efe[s]^{z\times n}$ such that $\deg(W(s))\leq d$ and
$\begin{bmatrix}P(s)\\W(s)\end{bmatrix}\approx Q(s)$.
\hfill $\Box$

\begin{remark}
  \label{remuvtocd}
  In Theorem \ref{theoprescr4} some of the conditions are expressed in terms of the row minimal indices. Occasionally, in the paper it is necessary to express them in terms of the column minimal indices. 

Assume that (\ref{eqinterfipolhom})-(\ref{eqdegsumpol}) hold. From (\ref{eqinterfipolhom}), (\ref{eqcmimajpol}), and (\ref{eqdegsumpol}) for $x=0$, we obtain 
$$
\sum_{i=1}^{m+z-r-x}v_i-\sum_{i=1}^{m-r}u_i+\sum_{i=1}^{r+x}\deg(\gamma_i)=
\sum_{i=1}^{n-r}c_i-\sum_{i=1}^{n-r-x}d_i+\sum_{i=1}^{r}\deg(\phi_i)+xd,
$$
therefore, from (\ref{eqdegsumpol}), and  taking into account (\ref{eqrmimajpol}) for $x=z$,   we get
\begin{equation}\label{eqdegsumpolequiv}
\begin{aligned}
  \sum_{i=1}^{r+x}\deg(\lcm(\phi_{i-x},\gamma_i))
  \leq \sum_{i=1}^{n-r}c_i-\sum_{i=1}^{n-r-x}d_i+\sum_{i=1}^{r}\deg(\phi_i)+xd,\\
  \mbox{with equality when $x=z$}.
  \end{aligned}
\end{equation}
Moreover,    (\ref{eqcmimajpol}) and  (\ref{eqrmimajpol}) hold for
   $\ba=(a_1, \dots, a_x)$ and $\bb=(b_1, \dots, b_{z-x})$ defined as 
\begin{equation}\label{eqdefabis}
\begin{array}{rl}
a_1=&
  \sum_{i=1}^{n-r}c_i-\sum_{i=1}^{n-r-x}d_i\\&+\sum_{i=1}^{r}\deg(\phi_i)
  -
\sum_{i=1}^{  r+x-1}\deg(\lcm(  \phi_{i-x+1},  \gamma_i))+(x-1)d,\\
a_j=&
\sum_{i=1}^{  r+x-j+1}\deg(\lcm(  \phi_{i-x+j-1},  \gamma_i))
-
\sum_{i=1}^{  r+x-j}\deg(\lcm(  \phi_{i-x+j},  \gamma_i))
-d,\\ & \hfill 2\leq j \leq x,
\end{array}
\end{equation}
\begin{equation}\label{eqdefbbis}
\begin{array}{rl}
 b_1=&
\sum_{i=1}^{n-r}c_i-\sum_{i=1}^{n-r-x}d_i+\sum_{i=1}^{r}\deg(\phi_i)-
\sum_{i=1}^{  r+x}\deg(\lcm(  \phi_{i-x-1},  \gamma_i))+xd,
\\
  b_j=&
\sum_{i=1}^{  r+x}\deg(\lcm(  \phi_{i-x-j+1},  \gamma_i)-
\sum_{i=1}^{  r+x}\deg(\lcm(  \phi_{i-x-j},  \gamma_i)),\\& \hfill 2\leq j \leq z-x.
\end{array}
\end{equation}

Conversely, (\ref{eqinterfipolhom})-(\ref{eqrmimajpol}) and (\ref{eqdegsumpolequiv}) 
with $\ba$ and $\bb$   defined as in 
(\ref{eqdefabis}) and (\ref{eqdefbbis}), respectively, imply (\ref{eqinterfipolhom})-(\ref{eqdegsumpol}),
 with $\ba$ and $\bb$   defined as in 
(\ref{eqdefa}) and (\ref{eqdefb}).
\end{remark}

\subsection{Prescription of finite and  infinite structures and column minimal indices}\label{subsec_fininfcol}

In this  subsection we solve  Problem \ref{problem} when the invariants to be achieved are the  invariant factors, the infinite elementary divisors and  the column minimal indices (equivalently, the homogeneous invariant factors  and the  column  minimal indices).
Previously, we need a  technical lemma.

Given two sequences of integers $\bu = (u_1, \dots, u_{p})$ and $\bb = (b_1, \dots, b_{y})$ the union, $\bu\cup \bb$,  is the decreasingly ordered sequence  of the $p+y$ integers of $\bu$ and $\bb$.

\begin{lemma}\label{lemmacup}
 Let $\bu = (u_1, \dots, u_{p})$ and $\bb = (b_1, \dots, b_{y})$    be sequences of  integers.  
Then
$$
\bu\cup \bb\prec'(\bu, \bb).
$$
  \end{lemma}
    \noindent
{\bf Proof}.
    Let $\bv=\bu\cup \bb$ and let
$k_j=\min \{i \; : \; b_j>u_{i}\}$, $1\leq j \leq y$.
We have $1\leq k_1\leq k_2\leq \dots\leq k_y\leq p+1$, and (recall that   if a condition is stated for $a\leq i\leq b$ with $a>b$, we understand that the condition disappears)
$$
\begin{array}{lll}
v_{i}=u_i\geq b_1, & 1\leq i \leq k_1-1, & 
v_{k_1}=b_1>u_{k_1},\\
v_{i}=u_{i-1}\geq b_2, & k_1+1\leq i \leq k_2, &
v_{k_2+1}=b_2>u_{k_2},\\
\vdots&\vdots & \vdots\\
v_{i}=u_{i-y+1}\geq b_y, & k_{y-1}+y-1\leq i \leq k_y+y-2, &
v_{k_y+y-1}=b_y>u_{k_y}, \\
v_{i}=u_{i-y}\geq b_{y+1}=-\infty, & k_{y}+y\leq i \leq p+y.& 
\end{array}
$$
In summary, putting $k_0=1$ and $k_{y+1}=p+1$, 
$$
\begin{array}{ll}
v_{i}=u_{i-j+1}\geq b_{j}, & k_{j-1}+j-1\leq i \leq k_{j}+j-2, \quad 1\leq j \leq y+1,  \\
v_{k_{j}+j-1}=b_{j}>u_{k_{j}}, & 1\leq j \leq y.\\
\end{array}
$$
It is clear that 
$
u_{i}\geq v_{i+y}$, $  1\leq i \leq p.
$

For $1\leq j\leq y$, let $h_j=\min\{i\; : \; u_{i-j+1}<v_i\}$. 
Let us see that $h_j=k_j+j-1$. For $i=k_j+j-1$
we obtain $u_{i-j+1}=u_{k_{j}}<b_j=v_{k_{j}+j-1}$; hence $h_j\leq k_j+j-1$. Moreover, if $i\leq k_j+j-2$, then
$i-j+1\leq  k_j-1$. Let $\ell=\min\{r\; :\; i-j+1\leq  k_r-1\}$. Then $\ell\leq j$,  
$k_{\ell-1}\leq i-j+1 \leq k_{\ell}-1$, and $u_{i-j+1}=v_{i+\ell-j}\geq v_i$; hence, $h_j=k_j+j-1$.
Furthermore, for $1\leq j\leq y$,
$$
\begin{array}{rl}
\sum_{i=1}^{h_j}v_i=&\sum_{i=1}^{k_j+j-1}v_i
=\sum_{r=1}^j\sum_{i=k_{r-1}+r-1}^{k_r+r-1}v_i\\=&\sum_{r=1}^jb_r+\sum_{r=1}^j\sum_{i=k_{r-1}+r-1}^{k_r+r-2}u_{i-r+1}
=\sum_{r=1}^jb_r+\sum_{r=1}^j\sum_{i=k_{r-1}}^{k_r-1}u_{i}\\=&\sum_{i=1}^jb_i+\sum_{i=1}^{k_j-1}u_{i}=\sum_{i=1}^jb_i+\sum_{i=1}^{h_j-j}u_{i}.
\end{array}
$$
Finally, notice that
$\sum_{i=1}^{p+y}v_i=\sum_{i=1}^{p}u_i+\sum_{i=1}^{y}b_i$.
\hfill $\Box$

\begin{theorem}
  \label{corprescrhifcmi}
Let $P(s)\in\efe[s]^{m\times n}$ be a polynomial matrix, $\deg(P(s))=d$, $\rank (P(s))=r$.
Let $\phi_1(s,t)\mid\cdots\mid\phi_r(s,t)$ be its
  homogeneous  invariant factors and 
  $\bc=(c_1,  \dots,  c_{n-r})$  its
column minimal indices.

Let $z$ and $x$ be integers such that $0\leq x\leq \min\{z, n-r\}$ and  
let  $\gamma_1(s, t)\mid\cdots\mid\gamma_{r+x}(s,t)$ be monic homogeneous polynomials 
and
$d_1 \geq  \dots  \geq  d_{n-r-x}$ non negative integers.
 There exists a polynomial matrix  $W(s)\in \efe[s]^{z\times n}$ such that $\deg(W(s))\leq d$, $\rank \left(\begin{bmatrix}P(s)\\W(s)\end{bmatrix}\right)=r+x$, and $\begin{bmatrix}P(s)\\W(s)\end{bmatrix}$ has
$\gamma_1(s, t)\mid\cdots\mid\gamma_{r+x}(s,t)$ as  homogeneous invariant factors and
$\bd=(d_1,\dots, d_{n-r-x})$  as column minimal indices
if and only if (\ref{eqinterfipolhom}),
 (\ref{eqcmimajpol}) and  (\ref{eqdegsumpolequiv})  hold, 
where $\ba = (a_1, \dots, a_x )$ is  defined as in (\ref{eqdefabis}).
\end{theorem}

\noindent
{\bf Proof}. 
The necessity follows directly  from Theorem \ref{theoprescr4} (see Remark \ref{remuvtocd}).
For the converse, let us assume that  (\ref{eqinterfipolhom}),
 (\ref{eqcmimajpol}) and (\ref{eqdegsumpolequiv}) hold.
Let $ \bu=(u_1, \dots,  u_{m-r})$  be the row minimal
indices of $P(s)$, where   $u_1 \geq \dots \geq u_{\eta}  >  u_{\eta +1}= \dots = u_{m-r}=0 $. Define
$\bb=( b_1, \dots, b_{z-x})$
as in (\ref{eqdefbbis})  and 
$\bv=(v_1, \dots, v_{m+z-r-x})=\bu \cup  \bb.$
Then, by  Lemma \ref{lemmacup}, condition (\ref{eqrmimajpol}) holds.
 
Denoting $\bar \eta=\#\{i\; : \; v_i>0\}$, clearly (\ref{eqetapol}) holds.  By Theorem \ref{theoprescr4} and  Remark \ref{remuvtocd}, the result follows. 
\hfill $\Box$

\subsection{Prescription of finite and  infinite structures and row minimal indices}\label{subsec_fininfrow}

 Now  we solve  Problem \ref{problem}  when the prescribed invariants are the homogeneous invariant factors  and the  row minimal indices.
We   use the following technical lemma.

\begin{lemma}
 \label{lemmaexistsd}
Let $\bc = (c_1, \dots, c_{q})$ and $\ba=(a_1, \dots, a_x)$   be sequences of  integers such that $q> x \geq 0$.
Let
$\ell=\min\{j\; : \; \sum_{i=1}^{j}c_i>\sum_{i=1}^{j}a_i\}$.
\begin{enumerate}
  \item\label{it1_lemmaexistsd}
If there exists a sequence of integers
$\bd=(d_1, \ldots, d_{q-x})$ such that  $\bc\prec'(\bd, \ba)$ then
\begin{equation}\label{eqsumca2}
\sum_{i=1}^{x+1}c_i-c_\ell\geq \sum_{i=1}^xa_i,
\end{equation}
\begin{equation}\label{eqsumca3}
\sum_{i=j+2}^{x+1}c_i\geq \sum_{i=j+1}^xa_i,\quad \ell\leq j \leq x-1.
\end{equation}
\item \label{it2_lemmaexistsd}
  If (\ref{eqsumca2}) and (\ref{eqsumca3}) hold, let
  \begin{equation}\label{defd}
d_1=\sum_{i=1}^{x+1}c_i-\sum_{i=1}^xa_i, \quad d_i=c_{i+x}, \quad 2\leq i \leq q-x.
\end{equation}
Then $d_1\geq  \dots \geq  d_{q-x}$  and $\bc\prec'(\bd, \ba)$, where $\bd=(d_1,  \dots, d_{q-x})$.
\end{enumerate}
  \end{lemma}

\begin{remark}\label{remell}
  Bearing in mind that $q\geq x+1$ and $a_{x+1}=-\infty$, we have  $\sum_{i=1}^{x+1}c_i>\sum_{i=1}^{x+1}a_i$; hence, $\ell$ is well defined, 
  $1\leq \ell\leq x+1$, and  $c_\ell \geq c_{x+1}$. If condition   (\ref{eqsumca2})  is satisfied, then  $\sum_{i=1}^xc_i\geq \sum_{i=1}^xa_i$.

   Moreover, if $\ell=x$ or $\ell=x+1$, then (\ref{eqsumca3}) vanishes, therefore it is trivially fulfilled.
\end{remark}
\noindent
{\bf Proof}.[Proof of Lemma \ref{lemmaexistsd}]\ 
 \begin{enumerate}
\item    Let us assume that there exists a sequence of integers
    $\bd=(d_1, \ldots, d_{q-x})$ such that  $\bc\prec'(\bd, \ba)$.
    It implies
    $\sum_{i=1}^qc_i= \sum_{i=1}^{q-x}d_i+ \sum_{i=1}^xa_i$ and $d_i\geq c_{i+x}$, $1\leq i\leq q-x$; hence 
    \begin{equation}\label{eqdleq}
\sum_{i=1}^jd_i= \sum_{i=1}^{q}c_i- \sum_{i=1}^xa_i-\sum_{i=j+1}^{q-x}d_i\leq \sum_{i=1}^{x+j}c_i- \sum_{i=1}^xa_i, \quad 0\leq j \leq q-x.
    \end{equation}
    If $\ell=x+1$, then (\ref{eqsumca3}) vanishes and, from (\ref{eqdleq}) for $j=0$, $\sum_{i=1}^{x+1}c_i-c_\ell= \sum_{i=1}^xc_i\geq \sum_{i=1}^xa_i$; i.e. 
(\ref{eqsumca2}) holds.

    If $\ell\leq x$,
    let $h_j=\min\{i\; : \; d_{i-j+1}<c_i\}$, $1\leq j\leq x$. We know that  $\sum_{i=1}^{h_j}c_i\leq \sum_{i=1}^{h_j-j}d_i+ \sum_{i=1}^ja_i$, and $j\leq h_j\leq q-x+j$,
    $1\leq j\leq x$ (see Remark \ref{rem_maj}.\ref{haches}). 
As 
$\sum_{i=1}^{\ell}c_i>\sum_{i=1}^{\ell}a_i$, it follows that $\ell<h_\ell$; hence $d_1\geq c_\ell$.
From (\ref{eqdleq}), $c_\ell\leq d_1\leq  \sum_{i=1}^{x+1}c_i- \sum_{i=1}^xa_i$; i.e., 
(\ref{eqsumca2}) also holds. Moreover, for $\ell\leq j \leq x -1$,  $d_1\geq c_\ell\geq c_j$; therefore $h_j>j$, and 
from (\ref{eqdleq})   
$$
\sum_{i=1}^{h_j}c_i\leq \sum_{i=1}^{h_j-j}d_i+ \sum_{i=1}^j a_i\leq \sum_{i=1}^{x+h_j-j}c_i- \sum_{i=1}^xa_i+ \sum_{i=1}^j a_i=\sum_{i=1}^{x+h_j-j}c_i- \sum_{i=j+1}^xa_i;
$$
hence
$
\sum_{i=1}^{x-j}c_{h_j+i}=\sum_{i=h_j+1}^{h_j+x-j}c_i\geq  \sum_{i=j+1}^xa_i.
$
As $h_j\geq j+1$, we have $c_{h_j+i}\leq c_{j+i+1}$, $1\leq i \leq x-j$. Thus, $\sum_{i=1}^{x-j}c_{h_j+i}\leq\sum_{i=1}^{x-j}c_{j+i+1}=\sum_{i=j+2}^{x+1}c_{i}$, from where we obtain
(\ref{eqsumca3}).

\item Let us assume that (\ref{eqsumca2}) and (\ref{eqsumca3}) hold and let $d_1,\ldots,d_{q-x}$ be defined as in (\ref{defd}). 
From (\ref{eqsumca2}), $d_1\geq c_\ell\geq c_{1+x}$; thus,
$d_1\geq  \dots \geq  d_{q-x}$ and $d_i\geq c_{i+x}$, $1\leq i\leq q-x$.
Moreover, $\sum_{i=1}^{q-x}d_i=\sum_{i=1}^{q}c_i-\sum_{i=1}^{x}a_i$.
Let $h_j=\min\{i\;: \; d_{i-j+1}<c_i\}$, $1\leq j\leq x$. It only remains to prove that, for $1\leq j \leq x$,
\begin{equation}\label{eqsumh}
\sum_{i=1}^{h_j}c_i\leq \sum_{i=1}^{h_j-j}d_i+ \sum_{i=1}^ja_i.
\end{equation}
We have $d_{i-j+1}\geq c_i$, $j\leq i \leq h_j-1$ and, since $h_j\geq j$, it follows that $c_{h_j} \leq c_j$.

Let $j\in\{1, \dots, \ell-1\}$.
By the definition of $\ell$,  $\sum_{i=1}^{j}c_i\leq\sum_{i=1}^{j}a_i$,  and 
$$\sum_{i=1}^{h_j}c_i=c_{h_j}+\sum_{i=1}^{j-1}c_i+\sum_{i=j}^{h_j-1}c_i\leq \sum_{i=1}^{j}c_i+\sum_{i=j}^{h_j-1}d_{i-j+1}\leq \sum_{i=1}^{j}a_i+\sum_{i=1}^{h_j-j}d_{i};
$$
i.e., (\ref{eqsumh}) holds for $1\leq j \leq \ell-1$.

Let $j\in\{\ell, \dots, x\}$.
Since $d_1\geq c_\ell\geq c_j$, we have $h_j>j$. Let $h_j-j=k$. Then $1\leq k\leq q-x$ and
$$
\begin{array}{ll}
& \sum_{i=1}^{h_j-j}d_i+ \sum_{i=1}^ja_i-\sum_{i=1}^{h_j}c_i \\
= &
\sum_{i=1}^{x+1}c_i- \sum_{i=1}^xa_i+\sum_{i=2}^{h_j-j}c_{i+x}+ \sum_{i=1}^ja_i-\sum_{i=1}^{h_j}c_i \\
= & 
\sum_{i=1}^{k+x}c_{i}-\sum_{i=1}^{k+j}c_i-\sum_{i=j+1}^xa_i=\sum_{i=k+j+1}^{k+x}c_{i}-\sum_{i=j+1}^xa_i.
\end{array}
$$
If $k=1$ then, from (\ref{eqsumca3}), we obtain $\sum_{i=k+j+1}^{k+x}c_{i}-\sum_{i=j+1}^xa_i\geq 0$; i.e. (\ref{eqsumh}) holds.
If $k>1$, then by the definition of $h_j$,  $d_{i+1}\geq c_{j+i}$ for $0\leq i \leq k-1$. Then
$c_{x+i+1}\geq c_{j+i}$ for $1\leq i \leq k-1$; hence
$c_{j+i}=c_{j+i+1}=\dots=c_{x+i+1}$ for $1\leq i \leq k-1$.
 For  $2\leq i \leq k-1$, we have $j+i-1\leq j+i\leq x+i$,  which means that these sequences overlap, therefore
$$
c_{j+1}=c_{j+2}=\dots =c_{j+k-1}=\dots=c_{x+k}.
$$
As a consequence, 
$\sum_{i=j+2}^{x+1}c_{i}=(x-j)c_{x+k}=\sum_{i=k+j+1}^{k+x}c_{i}$; hence
$$\sum_{i=1}^{h_j-j}d_i+ \sum_{i=1}^ja_i-\sum_{i=1}^{h_j}c_i=\sum_{i=j+2}^{x+1}c_{i}-\sum_{i=j+1}^xa_i.$$
From (\ref{eqsumca3}) we derive (\ref{eqsumh}).
\end{enumerate}
\hfill $\Box$

\begin{theorem}
  \label{corprescrhifrmi}
  Let $P(s)\in\efe[s]^{m\times n}$ be a polynomial matrix, $\deg(P(s))=d$, $\rank (P(s))=r$.
Let $\phi_1(s,t)\mid\cdots\mid\phi_r(s,t)$ be its
  homogeneous  invariant factors, 
  $\bc=(c_1,  \dots,  c_{n-r})$  its
column minimal indices, and $ \bu=(u_1, \dots,  u_{m-r})$  its row minimal
indices, where   $u_1 \geq \dots \geq u_{\eta}  >  u_{\eta +1}= \dots = u_{m-r}=0 $.

Let $z$ and $x$ be integers such that $0\leq x\leq \min\{z, n-r\}$ and 
 let  $\gamma_1(s, t)\mid\cdots\mid\gamma_{r+x}(s,t)$ be  monic homogeneous polynomials 
 and
 $v_1 \geq  \dots  \geq  v_{\bar \eta} > v_{\bar \eta+1}=\dots= v_{m+z-r-x}=0$ be non negative integers.
 Let
$\ba = (a_1, \dots, a_x )$ and $\bb = (b_1, \dots,  b_{z-x} )$
be as in (\ref{eqdefa}) and (\ref{eqdefb}), respectively.
\begin{enumerate}
\item
  If $x=n-r$, 
 there exists a polynomial matrix  $W(s)\in \efe[s]^{z\times n}$ such that $\deg(W(s))\leq d$, $\rank\left(\begin{bmatrix}P(s)\\W(s)\end{bmatrix}\right)=r+x$,
 and 
 $\begin{bmatrix}P(s)\\W(s)\end{bmatrix}$ has 
$\gamma_1(s, t)\mid\cdots\mid\gamma_{r+x}(s,t)$ as homogeneous invariant factors and
$\bv=(v_1,\dots, v_{m+z-r-x})$ as row minimal indices
if and only if (\ref{eqinterfipolhom}),
(\ref{eqetapol}), (\ref{eqrmimajpol}), (\ref{eqdegsumpol})
 and 
\begin{equation}\label{eqcpreca}\bc \prec \ba.\end{equation}
\item
  If $x<n-r$, 
 there exists a polynomial matrix  $W(s)\in \efe[s]^{z\times n}$ such that $\deg(W(s))\leq d$, $\rank\left(\begin{bmatrix}P(s)\\W(s)\end{bmatrix}\right)=r+x$, 
 and  $\begin{bmatrix}P(s)\\W(s)\end{bmatrix}$ has 
$\gamma_1(s, t)\mid\cdots\mid\gamma_{r+x}(s,t)$ as homogeneous invariant factors and
$\bv=(v_1,\dots, v_{m+z-r-x})$ as row minimal indices
if and only if (\ref{eqinterfipolhom}),
(\ref{eqetapol}), (\ref{eqrmimajpol}), (\ref{eqdegsumpol}),  (\ref{eqsumca2}) and (\ref{eqsumca3}) hold.
\end{enumerate}

\end{theorem}

\noindent
{\bf Proof}.\
\begin{enumerate}
\item
  Assume that $x=n-r$. The result follows from Theorem \ref{theoprescr4}
taking  into account that if $\bd=\emptyset$ then (\ref{eqcmimajpol}) and (\ref{eqcpreca}) are equivalent (see Remark \ref{rem_maj}.\ref{gm_pc}).
  
\item Assume that $x<n-r$.
  The necessity follows directly  from Theorem \ref{theoprescr4} and Lemma \ref{lemmaexistsd}.\ref{it1_lemmaexistsd}.
Conversely, let us assume that 
(\ref{eqinterfipolhom}),
(\ref{eqetapol}), \ref{eqrmimajpol}), (\ref{eqdegsumpol}), (\ref{eqsumca2}) and (\ref{eqsumca3}) hold.
From  (\ref{eqsumca2}) and (\ref{eqsumca3}), 
by Lemma \ref{lemmaexistsd}.\ref{it2_lemmaexistsd}
there exists a  sequence of integers
$\bd=(d_1, \ldots, d_{n-r-x})$   satisfying (\ref{eqcmimajpol}).
We have $d_{n-r-x}\geq c_{n-r}\geq 0$. 
By Theorem \ref{theoprescr4},  the result follows.
\end{enumerate}
\hfill $\Box$

\subsection{Prescription of finite   and/or infinite structures}\label{subsec_fininf}

This  subsection is devoted to solve   Problem \ref{problem}  when the  invariant factors and/or the  infinite elementary divisors are prescribed. In Theorem \ref{theoprescrhom}, both the finite and infinite structures are prescribed, in Theorem \ref{theoprescrfif} we only prescribe the finite  structure,  and  in Theorem \ref{theoprescried} we only prescribe the infinite structure.
We also start with a technical lemma.

\begin{lemma}\label{lemmagmm}
Let $\bc = (c_1, \dots, c_{q})$,  $\bd = (d_1, \dots, d_{q-x})$ and $\ba=(a_1, \dots, a_x)$   be sequences of  integers. 
If
$\bc \prec' (\bd, \ba)$ then $(c_1, \dots, c_{x})\prec (a_1+K, a_2, \dots, a_x)$, where $K=\sum_{i=1}^{q-x}(d_i-c_{i+x})$.
  \end{lemma}
\noindent
{\bf Proof}. For $1\leq j\leq x$, let $h_j=\min\{i\; : \; d_{i-j+1}<c_i\}$. We know that $j \leq h_j\leq q-x+j$ (see Remark \ref{rem_maj}.\ref{haches}) and 
$
\sum_{i=1}^{h_j}c_i\leq \sum_{i=1}^{h_j-j}d_i+ \sum_{i=1}^{j}a_i;
$
hence
$$
\begin{array}{rl}
\sum_{i=1}^{j}c_i\leq& \sum_{i=1}^{j}a_i+ \sum_{i=1}^{h_j-j}d_i-\sum_{i=j+1}^{h_j}c_i=
\sum_{i=1}^{j}a_i+ \sum_{i=1}^{h_j-j}d_i-\sum_{i=1}^{h_j-j}c_{i+j}\\\leq& \sum_{i=1}^{j}a_i+ \sum_{i=1}^{h_j-j}d_i-\sum_{i=1}^{h_j-j}c_{i+x}, \quad 1\leq j \leq x.
\end{array}
$$
Moreover, as $d_i\geq c_{i+x}$ for $1\leq i \leq q-x$, we conclude that $\sum_{i=1}^{h_j-j}(d_i-c_{i+x})\leq\sum_{i=1}^{q-x}(d_i-c_{i+x})=K$. Therefore, we obtain
$$
\sum_{i=1}^{j}c_i\leq \sum_{i=1}^{j}a_i+K, \quad 1\leq j \leq x.
$$
Finally, notice that 
$$
\sum_{i=1}^{x}c_i=\sum_{i=1}^{x}a_i+\sum_{i=1}^{q-x}d_i-\sum_{i=1}^{q-x}c_{i+x}=\sum_{i=1}^{x}a_i+K,
$$
is also satisfied.
\hfill $\Box$

\begin{theorem}
  \label{theoprescrhom}
  Let $P(s)\in\efe[s]^{m\times n}$, $\deg(P(s))=d$, $\rank (P(s))=r$.
Let $\phi_1(s,t)\mid\cdots\mid\phi_r(s,t)$ be its
  homogeneous  invariant factors, 
  $\bc=(c_1,  \dots,  c_{n-r})$  its
column minimal indices, and $ \bu=(u_1, \dots,  u_{m-r})$  its row minimal
indices, where   $u_1 \geq \dots \geq u_{\eta}  >  u_{\eta +1}= \dots = u_{m-r}=0 $.

Let $z$ and $x$ be integers such that $0\leq x\leq \min\{z, n-r\}$ and  
let  $\gamma_1(s, t)\mid\cdots\mid\gamma_{r+x}(s,t)$ be  monic homogeneous polynomials.
\begin{enumerate}
 \item  \label{ittheoprescrhom1} If  $x<z$ or $x=z=n-r$, then 
there exists a polynomial matrix  $W(s)\in \efe[s]^{z\times n}$ such that  $\deg(W(s))\leq d$, $\rank\left(\begin{bmatrix}P(s)\\W(s)\end{bmatrix}\right)=r+x$
 and
$\begin{bmatrix}P(s)\\W(s)\end{bmatrix}$ has 
$\gamma_1(s, t)\mid\cdots\mid\gamma_{r+x}(s,t)$ as homogeneous invariant factors
 if and only if  (\ref{eqinterfipolhom}) and
 \begin{equation}\label{eqfromhifj}
 \begin{aligned}
   \begin{array}{l}
 \sum_{i=1}^{r+x-j}\deg(\lcm(\phi_{i-x+j},\gamma_i))+\sum_{i=1}^{m-r}u_i+\sum_{i=1}^{j}c_i+\sum_{i=x+1}^{n-r}c_i
 \\ \leq (r+x-j)d, \quad  0\leq j \leq x-1,
  \end{array}
  \\
    \mbox{with equality for $j=0$ when $x=z=n-r$}.
 \end{aligned}
 \end{equation}

\item \label{ittheoprescrhom2}
  If  $x=z<n-r$, then  there exists a polynomial matrix  $W(s)\in \efe[s]^{z\times n}$ such that  $\deg(W(s))\leq d$, $\rank\left(\begin{bmatrix}P(s)\\W(s)\end{bmatrix}\right)=r+x$
 and
$\begin{bmatrix}P(s)\\W(s)\end{bmatrix}$ has 
$\gamma_1(s, t)\mid\cdots\mid\gamma_{r+x}(s,t)$ as homogeneous invariant factors
 if and only if  (\ref{eqinterfipolhom}),
  \begin{equation}\label{eqnew1}
\sum_{i=1}^{x+1}c_i-c_\ell\geq \sum_{i=1}^{r+x}\deg(\gamma_i)-\sum_{i=1}^{r}\deg(\phi_i)-xd,
  \end{equation}
  \begin{equation}\label{eqnew2}
    \begin{array}{l}
 \sum_{i=j+2}^{x+1}c_i\geq\sum_{i=1}^{r+x-j}\deg(\lcm(\phi_{i-x+j},\gamma_i))-\sum_{i=1}^{r}\deg(\phi_i)-(x-j)d,\\
 \hfill \ell\leq j \leq x-1,
\end{array}
\end{equation}
where $\ell=\min\{j\, :\,\sum_{i=1}^{j}c_i>\sum_{i=1}^{r+x}\deg(\gamma_i)-\sum_{i=1}^{r+x-j}\deg(\lcm(\phi_{i-x+j},\gamma_i))-jd \}$.
\end{enumerate}
\end{theorem}

\noindent
{\bf Proof}.\, 
  \begin{enumerate}
  \item   Case  $x<z$ or $x=z=n-r$.
Assume that 
there exists  $W(s)\in \efe[s]^{z\times n}$ such that $\deg(W(s))\leq d$, $\rank\left(\begin{bmatrix}P(s)\\W(s)\end{bmatrix}\right)=r+x$ 
and
 $\begin{bmatrix}P(s)\\W(s)\end{bmatrix}$ has 
$\gamma_1(s, t)\mid\cdots\mid\gamma_{r+x}(s,t)$ as homogeneous invariant factors.
Let 
$\bd=(d_1,\dots, d_{n-r-x})$ be the 
column minimal indices of $\begin{bmatrix}P(s)\\W(s)\end{bmatrix}$.
By Theorem \ref{corprescrhifcmi}, (\ref{eqinterfipolhom}),
 (\ref{eqcmimajpol}) and (\ref{eqdegsumpolequiv}) hold, 
where $\ba = (a_1, \dots, a_x )$ is  defined as in (\ref{eqdefabis}).
From (\ref{eqdegsumpolequiv}) 
we obtain
$$
  \sum_{i=1}^{r+x}\deg(\lcm(\phi_{i-x},\gamma_i))
  +\sum_{i=1}^{m-r}u_i+\sum_{i=1}^{n-r-x}d_i
  \leq (r+x)d,$$
with equality if $x=z$.
From (\ref{eqcmimajpol})  we get $\sum_{i=1}^{n-r-x}d_i\geq \sum_{i=1}^{n-r-x}c_{i+x}=\sum_{i=x+1}^{n-r}c_i$.
Therefore,  (\ref{eqfromhifj}) 
holds for $j=0$.

For $1\leq j\leq x -1$, from (\ref{eqcmimajpol}), Lemma \ref{lemmagmm}, and the definition of $a_j$,  we obtain
$$
\begin{array}{rl}
\sum_{i=1}^{j}c_i\leq & \sum_{i=1}^{j}a_i+\sum_{i=1}^{n-r-x}d_i-\sum_{i=x+1}^{n-r}c_{i}\\
=& \sum_{i=1}^{n-r}c_i-\sum_{i=1}^{n-r-x} d_i+\sum_{i=1}^{r}\deg( \phi_i)+xd\\
& -
\sum_{i=1}^{  r+x-j}\deg(\lcm(  \phi_{i-x+j},  \gamma_i))-jd+\sum_{i=1}^{n-r-x}d_i\\&-\sum_{i=x+1}^{n-r}c_{i}\\=&
(r+x-j)d-\sum_{i=1}^{m-r}u_i-
\sum_{i=1}^{  r+x-j}\deg(\lcm(  \phi_{i-x+j},  \gamma_i))
\\&-\sum_{i=x+1}^{n-r}c_{i}.
\end{array}
$$
Thus, (\ref{eqfromhifj}) holds.

Conversely, 
assume that (\ref{eqinterfipolhom}) and (\ref{eqfromhifj}) hold.
Define
$$
\begin{array}{rl}
\hat a_1=&
 \sum_{i=1}^{x}c_i +\sum_{i=1}^{r}\deg(\phi_i)
 -
 \sum_{i=1}^{  r+x-1}\deg(\lcm(  \phi_{i-x+1},  \gamma_i))+(x-1)d, \\
 a_j=&
\sum_{i=1}^{  r+x-j+1}\deg(\lcm(  \phi_{i-x+j-1},  \gamma_i))
-
\sum_{i=1}^{  r+x-j}\deg(\lcm(  \phi_{i-x+j},  \gamma_i))
\\&-d, \hfill 2\leq j \leq x.\end{array}
$$
By condition (\ref{eqfromhifj}) for $j=0$,
$$(r+x)d \geq \sum_{i=1}^{r+x}\deg(\lcm(\phi_{i-x},\gamma_i))+rd-\sum_{i=1}^{r}\deg(\phi_i)-\sum_{i=1}^{x}c_i 
;$$
hence
$$\hat a_1\geq 
\sum_{i=1}^{r+x}\deg(\lcm(\phi_{i-x},\gamma_i))-\sum_{i=1}^{r+x-1}\deg(\lcm(\phi_{i-x+1},\gamma_i))-d.$$
 By Remark \ref{remdecreasing}.\ref{remdecreasing1},  we have  
$\hat a_1\geq a_2\geq\dots \geq a_{x}$. 
Let  $\hat \ba=(\hat a_1,  a_2,\dots, a_{x})$, then
for $1\leq j \leq x$, 
$$
\begin{array}{rl}
&\hat a_1+\sum_{i=2}^{j}a_i\\ =&\sum_{i=1}^{x}c_i +\sum_{i=1}^{r}\deg(\phi_i)
 -\sum_{i=1}^{  r+x-j}\deg(\lcm(  \phi_{i-x+j},  \gamma_i))+(x-j)d\\=&
(r+x-j)d-\sum_{i=x+1}^{n-r}c_i-\sum_{i=1}^{m-r}u_i-
 \sum_{i=1}^{  r+x-j}\deg(\lcm(  \phi_{i-x+j},  \gamma_i)).
 \end{array}
$$
From (\ref{eqfromhifj})  we obtain 
$$\hat a_1+\sum_{i=2}^{j}a_i\geq \sum_{i=1}^{ j}   c_i,\quad 1\leq j \leq x-1.$$
Moreover, 
$\hat a_1+\sum_{i=2}^{ x}a_i=\sum_{i=1}^{ x}c_i$.

If $n=r+x$ then $\bc \prec \hat \ba$, and let $\bd=\emptyset$ so that $\bc\prec'(\bd, \hat  \ba)$ holds. Otherwise, if $n>r+x$ by  Lemma \ref{lemmaexistsd} 
 there exists
a sequence of integers
$\bd=(d_1, \ldots, d_{n-r-x})$ such that $\bc\prec'(\bd, \hat  \ba)$, 
$d_i= c_{i+x}$ for $2\leq i \leq n-r-x$, and $d_1=\sum_{i=1}^{x+1}c_i-\hat{a}_1-\sum_{i=2}^{x}a_i=c_{x+1}$.

Let  $\ba=(a_1, \dots, a_{x})$ 
be defined as in (\ref{eqdefabis}). 
 Then $a_1=\hat a_1$,  therefore  (\ref{eqcmimajpol}) holds.
 From (\ref{eqfromhifj}) for $j=0$ we obtain
 $$
\begin{array}{rl}
 \sum_{i=1}^{r+x}\deg(\lcm(\phi_{i-x},\gamma_i))\hspace{-1pt}
 \leq\hspace{-2pt} & 
 \sum_{i=1}^{r}\deg(\phi_i)+\sum_{i=1}^{n-r}c_i-\sum_{i=1}^{n-r-x}c_{i+x} +xd\\\hspace{-1pt}=\hspace{-2pt}&\sum_{i=1}^{r}\deg(\phi_i)+\sum_{i=1}^{n-r}c_i-\sum_{i=1}^{n-r-x}d_{i} +xd,\end{array}
$$
with equality if $x=z=n-r$, i.e.,  (\ref{eqdegsumpolequiv}) is satisfied in this case.
By Theorem \ref{corprescrhifcmi},  the result follows.
\item
  Case $x=z<n-r$.   As $x=z$,  observe that if there exits $W(s)\in \efe[s]^{z\times n}$ such that $\rank\left(\begin{bmatrix}P(s)\\W(s)\end{bmatrix}\right)=r+x$,  then  the row minimal indices of 
$\begin{bmatrix}P(s)\\W(s)\end{bmatrix}$  are the row minimal indices of $P(s)$, i.e.,  $\bv=\bu$. For the sufficiency we prescribe $\bv=\bu$.  The result follows from Theorem  \ref{corprescrhifrmi}.
\end{enumerate}
\hfill $\Box$

\begin{remark}\label{4.1819-4.17}
We would like to remark that if $x=z<n-r$, conditions (\ref{eqinterfipolhom}), (\ref{eqnew1}) and (\ref{eqnew2}) imply (\ref{eqfromhifj}). 
Recall that $1 \leq \ell \leq x+1$ (see Remark \ref{remell}).

For $0\leq j<\ell$, from the definition of $\ell$ and (\ref{eqnew1}), 
$$\begin{array}{l}\sum_{i=1}^{r+x-j}\deg(\lcm(\phi_{i-x+j},\gamma_i))+\sum_{i=1}^{m-r} u_i+\sum_{i=1}^{j}c_i+\sum_{i=x+1}^{n-r} c_i\\\leq 
\sum_{i=1}^{r+x}\deg(\gamma_i)+\sum_{i=1}^{m-r}u_i+\sum_{i=x+1}^{n-r} c_i-jd\\
\leq
(r+x-j)d-c_{\ell}+c_{x+1}\leq (r+x-j)d.\end{array}$$
For $\ell \leq j\leq x-1$, from (\ref{eqnew2}),
$$\begin{array}{l} \sum_{i=1}^{r+x-j}\deg(\lcm(\phi_{i-x+j},\gamma_i))+\sum_{i=1}^{m-r} u_i+\sum_{i=1}^{j}c_i+\sum_{i=x+1}^{n-r} c_i\\ \leq
(r+x-j)d-c_{j+1}+c_{x+1}\leq (r+x-j)d.\end{array}$$
\end{remark}

When we only prescribe the invariant factors, we obtain the next theorem.
The solution to the row completion is now characterized in terms of two conditions: the interlacing inequalities (\ref{eqinterfipol}), and an  additional condition (\ref{eqfromifj}). In fact, condition
(\ref{eqinterfipol}) can be derived from \cite{Sa79, Th79}, where the problem of prescription of the invariant factors of a polynomial matrix when a submatrix is fixed was solved. The appearance of condition (\ref{eqfromifj})  is due to  degree condition on the completion matrix $W(s)$, restriction that was not taken into account in \cite{Sa79, Th79}.

\begin{theorem}\label{theoprescrfif}
  Let $P(s)\in\efe[s]^{m\times n}$,  $\deg(P(s))=d$, $\rank (P(s))=r$.
  Let $\alpha_1(s)\mid\cdots\mid\alpha_r(s)$ be its invariant factors,
$e_1\leq\cdots\leq e_r$
its partial multiplicities of $\infty$, 
$\bc=(c_1,  \dots,  c_{n-r})$  its
column minimal indices, and $ \bu=(u_1, \dots,  u_{m-r})$  its row minimal
indices.

Let $z$ and $x$ be integers such that $0\leq x\leq \min\{z, n-r\}$ and 
 let  $\beta_1(s)\mid\cdots\mid\beta_{r+x}(s)$ be monic polynomials.
There exists a polynomial matrix  $W(s)\in \efe[s]^{z\times n}$ such that  $\deg(W(s))\leq d$, $\rank\left(\begin{bmatrix}P(s)\\W(s)\end{bmatrix}\right)=r+x$ 
and $\begin{bmatrix}P(s)\\W(s)\end{bmatrix}$ has
$\beta_1(s)\mid\cdots\mid\beta_{r+x}(s)$ as invariant factors  
 if and only if
\begin{equation}\label{eqinterfipol}
 \beta_i(s)\mid \alpha_i(s)\mid \beta_{i+z}(s),\quad 1\leq i \leq r,
\end{equation}
\begin{equation}\label{eqfromifj}
  \begin{array}{l}
 \sum_{i=1}^{r+x-j}\deg(\lcm(\alpha_{i-x+j},\beta_i))+\sum_{i=1}^{r}e_i+\sum_{i=1}^{m-r}u_i+\sum_{i=1}^{j}c_i+\sum_{i=x+1}^{n-r}c_i
 \\ \leq (r+x-j)d, \quad  0\leq j \leq x-1.
  \end{array}
\end{equation}
\end{theorem}

\noindent
{\bf Proof}.
Let $\phi_1(s,t)\mid \dots \mid \phi_r(s,t)$ be the homogeneous invariant factors of $P(s)$.
Assume that 
there exists  $W(s)\in \efe[s]^{z\times n}$ such that  $\deg(W(s))\leq d$, $\rank\left(\begin{bmatrix}P(s)\\W(s)\end{bmatrix}\right)=r+x$, 
and $\begin{bmatrix}P(s)\\W(s)\end{bmatrix}$ has
$\beta_1(s)\mid\cdots\mid\beta_{r+x}(s)$ as  invariant factors.
Let $f_1\le \dots \le f_{r+x}$ be the partial multiplicities of $\infty$
and 
  $\gamma_1(s,t)\mid \dots \mid \gamma_{r+x}(s,t)$ the homogeneous invariant factors of
  $\begin{bmatrix}P(s)\\W(s)\end{bmatrix}$. 
By Theorem \ref{theoprescrhom} and Remark \ref{4.1819-4.17}, (\ref{eqinterfipolhom}) and 
(\ref{eqfromhifj}) hold.
Condition (\ref{eqinterfipolhom}) implies (\ref{eqinterfipol}).  
For $0\leq j \leq x-1$ we have
$$
\begin{array}{rl}
\sum_{i=1}^{r+x-j}\deg(\lcm(\phi_{i-x+j},\gamma_i))
\geq & \sum_{i=1}^{r+x-j}\deg(\lcm(\alpha_{i-x+j},\beta_i))+\sum_{i=1}^{r+x-j}e_{i-x+j}\\
= &\sum_{i=1}^{r+x-j}\deg(\lcm(\alpha_{i-x+j},\beta_i))+\sum_{i=1}^{r}e_{i};\\
\end{array}
$$
hence (\ref{eqfromhifj}) implies (\ref{eqfromifj}).

Conversely, 
assume that (\ref{eqinterfipol}) and (\ref{eqfromifj}) hold.
  Let $$q=(r+x)d-\sum_{i=1}^{r+x}\deg(\lcm(\alpha_{i-x},\beta_i))-\sum_{i=1}^{r}e_i-\sum_{i=1}^{m-r}u_i-\sum_{i=x+1}^{n-r}c_i.$$
  If $x=0$ then $q=0$ and if $x>0$, from (\ref{eqfromifj})  we obtain $q\geq 0$.
Define
$$
\begin{array}{l}
f_i=0, \quad  1\leq i \leq x,\\
f_{i+x}=e_{i}, \quad 1\leq i \leq r-1,\\f_{r+x}=e_{r}+q,\\ 
\end{array}
$$
and
$$\gamma_i(s,t)=t^{ f_i}t^{\deg(\beta_i)} \beta_i\left(\frac{s}{t}\right), \quad 1\leq i \leq  r+x.$$
We have $f_1\leq \dots\leq f_{r+x}$, and 
$
f_i\leq  e_i\leq f_{i+z}, 1\leq i \leq  r.
$
Thus, $\gamma_1(s,t)\mid \dots \mid \gamma_{r+x}(s,t)$ and from
(\ref{eqinterfipol}) we derive (\ref{eqinterfipolhom}).
If $x>1$, 
for $1\leq j \leq x$, we have $f_{i+x-j}\leq f_{i+x-1}=e_{i-1}\leq e_i$, $1\leq i \leq r$, therefore $
\sum_{i=1}^{r+x-j}\max\{e_{i-x+j}, f_i\}=\sum_{i=1}^{r}e_i
$, $1\leq j  \leq r$.
For $j=0$ we obtain $\sum_{i=1}^{r+x}\max\{e_{i-x}, f_i\}=\sum_{i=1}^{r}e_i+q$.
Thus,
$$
\begin{array}{rl}
  \sum_{i=1}^{r+x}\deg(\lcm(\phi_{i-x},\gamma_i))=&\sum_{i=1}^{r+x}\deg(\lcm(\alpha_{i-x},\beta_i))+\sum_{i=1}^{r}e_i+q\\=&
(r+x)d-\sum_{i=1}^{m-r}u_i-\sum_{i=x+1}^{n-r}c_i,\\
\sum_{i=1}^{r+x-j}\deg(\lcm(\phi_{i-x+j},\gamma_i))=&\sum_{i=1}^{r+x-j}\deg(\lcm(\alpha_{i-x+j},\beta_i))+\sum_{i=1}^{r}e_i,\\& \hfill 1 \leq j \leq x-1,
\end{array}
$$
hence  from (\ref{eqfromifj})
we obtain
(\ref{eqfromhifj}).

If $x<z$ or $x=z=n-r$ the result follows from  Theorem \ref{theoprescrhom} (item \ref{ittheoprescrhom1}).

If $x=z<n-r$, from (\ref{eqinterfipolhom}) we have
$$\sum_{i=1}^{r+x}\deg(\gamma_i)=\sum_{i=1}^{r+x}\deg(\lcm(\phi_{i-x},\gamma_i))=(r+x)d-\sum_{i=1}^{m-r}u_i-\sum_{i=x+1}^{n-r}c_i.$$ Therefore,
$$\sum_{i=1}^{r+x}\deg(\gamma_i)-\sum_{i=1}^{r}\deg(\lcm(\phi_{i}, \gamma_i))-xd=\sum_{i=1}^{x}c_i,$$ and 
from (\ref{eqfromhifj}) we obtain
$$
\begin{array}{rl}
&\sum_{i=1}^{r+x}\deg(\gamma_i)-\sum_{i=1}^{r+x-j}\deg(\lcm(\phi_{i-x+j},\gamma_i))-jd\\ 
  &=(r+x-j)d-\sum_{i=1}^{m-r}u_i-\sum_{i=x+1}^{n-r}c_i-\sum_{i=1}^{r+x-j}\deg(\lcm(\phi_{i-x+j},\gamma_i))\\
  & \geq \sum_{i=1}^jc_i, \hfill 1\leq j \leq x-1.
\end{array}
$$
Thus, (\ref{eqnew1}) holds and  (\ref{eqnew2}) vanishes. The result follows from 
 Theorem \ref{theoprescrhom} (item \ref{ittheoprescrhom2}).
\hfill $\Box$

If we only prescribe the infinite elementary divisors we obtain the following theorem.

\begin{theorem}\label{theoprescried}
  Let $P(s)\in\efe[s]^{m\times n}$, $\deg(P(s))=d$, $\rank (P(s))=r$.
  Let $\alpha_1(s)\mid\cdots\mid\alpha_r(s)$ be its invariant factors,
$e_1\le \dots \le e_{r}$ its partial multiplicities of $\infty$,
$\bc=(c_1,  \dots,  c_{n-r})$  its
column minimal indices, and $ \bu=(u_1, \dots,  u_{m-r})$  its row minimal
indices.

Let $z$ and $x$ be integers such that $0\leq x\leq \min\{z, n-r\}$ and 
 let   $f_1\leq \cdots\leq f_{r+x}$ be non negative integers.
There exists a polynomial matrix  $W(s)\in \efe[s]^{z\times n}$ such that  $\deg(W(s))\leq d$, $\rank\left(\begin{bmatrix}P(s)\\W(s)\end{bmatrix}\right)=r+x$
 and $\begin{bmatrix}P(s)\\W(s)\end{bmatrix}$ has
$f_1\le \dots \le f_{r+x}$ as partial multiplicities of $\infty$
 if and only if 
\begin{equation*}
f_i\leq  e_i\leq f_{i+z},\quad 1\leq i \leq  r,    
\end{equation*}
\begin{equation*}\label{eqfromidej}
  \begin{array}{l}
 \sum_{i=1}^{r+x-j}\max\{e_{i-x+j},f_i\}+\sum_{i=1}^{r}\deg(\alpha_i)+\sum_{i=1}^{m-r}u_i+\sum_{i=1}^{j}c_i+\sum_{i=x+1}^{n-r}c_i
 \\ \leq (r+x-j)d, \quad  0\leq j \leq x-1.
  \end{array}
\end{equation*}
\end{theorem}
\noindent
{\bf Proof}.
The proof is  analogous to that of the Theorem 
 \ref{theoprescrfif}.
To prove the sufficiency, define
$
  \beta_i(s)=1, 
  1\leq i \leq x$, 
  $\beta_{i+x}(s)=\alpha_{i}(s)$, 
  $1\leq i \leq r-1$ and $\beta_{r+x}(s)=\alpha_{r}(s)\tau(s)$, where 
$\tau(s)$ is a monic polynomial of 
$\deg(\tau)=(r+x)d-\sum_{i=1}^{r}\deg(\alpha_i)-\sum_{i=1}^{r+x}\max\{e_{i-x},f_i\}-\sum_{i=1}^{m-r}u_i-\sum_{i=x+1}^{n-r}c_i.$
\hfill $\Box$

\section{Conclusions and future work}\label{secconclusions}
In this work we have solved the problem of row (column) completion of a polynomial matrix of bounded degree when the eigenstructure has been prescribed (Subsection \ref{subsec_eigen}). We also have presented the  particular cases  where  we prescribe the finite and infinite structures and the column (row) minimal indices (Subsections \ref{subsec_fininfcol} and \ref{subsec_fininfrow}), and the finite and/or infinite structures (Subsection \ref{subsec_fininf}). 
In a future paper we will accomplish the study of the remaining cases, i.e.,  the prescription of the column and/or row minimal indices,  of both the column and row minimal indices and the finite (respectively, infinite) structure,  of the column minimal indices and the finite (respectively, infinite) structure and, finally,  of  the row minimal indices and the finite (respectively, infinite) structure.

\bibliographystyle{acm}
\bibliography{references}

\begin{thebibliography}{10}

\bibitem{Ba89}
{\sc {Baraga\~na}, I.}
\newblock Interlacing inequalities for regular pencils.
\newblock {\em Linear Algebra Appl. 121\/} (1989), 521--531.

\bibitem{BaZa90}
{\sc {Baraga\~na}, I., and Zaballa, I.}
\newblock Column completion of a pair of matrices.
\newblock {\em Linear Multilinear Algebra 27\/} (1990), 243--273.

\bibitem{Bhms11}
{\sc Betcke, T., Higham, N.~J., Mehrmann, V., {Schr\"oder}, C., and Tisseur,
  F.}
\newblock Nlevp: A collection of nonlinear eigenvalue problems.
\newblock {\em ACM Transactions on Mathematical Software 39}, 2 (2013), 1--28.

\bibitem{BoDo94}
{\sc Boley, D.~L., and Dooren, P.~V.}
\newblock Placing zeroes and the {Kronecker} canonical form.
\newblock {\em Circuits Systems and Signal Process 13\/} (1994), 783--802.

\bibitem{CaSi91}
{\sc Cabral, I., and Silva, F.~C.}
\newblock Unified theorems on completions of matrix pencils.
\newblock {\em Linear Algebra Appl. 159\/} (1991), 43--54.

\bibitem{CaSi92}
{\sc Cabral, I., and Silva, F.~C.}
\newblock Similarity invariants of completions of submatrices.
\newblock {\em Linear Algebra Appl. 169\/} (1992), 151--161.

\bibitem{Di74}
{\sc da~Silva, J. A.~D.}
\newblock Matrices with prescribed entries and characteristic polynomial.
\newblock {\em Proc. Amer Math. Soc. 45}, 1 (July 1974), 31--37.

\bibitem{Ol69}
{\sc de~Oliveira, G.}
\newblock Matrices with prescribed characteristic polynomial and a prescribed
  submatrix. {I}.
\newblock {\em Pacific Journal of Mathematics 29}, 3 (1969), 653--661.

\bibitem{DeDoVa15}
{\sc {De Ter\'an}, F., Dopico, F.~M., and Dooren, P.~V.}
\newblock Matrix polynomials with completely prescribed eigenstructure.
\newblock {\em SIAM J. Matrix Anal. Appl. 36}, 1 (2015), 302--328.

\bibitem{DeDoMa14}
{\sc {De Ter\'an}, F., Dopico, F.~M., and Mackey, D.~S.}
\newblock Spectral equivalence of matrix polynomials and the index sum theorem.
\newblock {\em Linear Algebra Appl. 459\/} (2014), 264--333.

\bibitem{Do05}
{\sc Dodig, M.}
\newblock Feedback invariants of matrices with prescribed rows.
\newblock {\em Linear Algebra Appl. 405\/} (2005), 121--154.

\bibitem{Do08}
{\sc Dodig, M.}
\newblock Matrix pencils completions problems.
\newblock {\em Linear Algebra Appl. 428\/} (2008), 259--304.

\bibitem{Do10}
{\sc Dodig, M.}
\newblock Explicit solution of the row completion problem for matrix pencils.
\newblock {\em Linear Algebra Appl. 432\/} (2010), 1299--1309.

\bibitem{Do13}
{\sc Dodig, M.}
\newblock Completion up to a matrix pencil with column minimal indices as the
  only nontrivial kronecker invariants.
\newblock {\em Linear Algebra Appl. 438\/} (2013), 3155--3173.

\bibitem{Do22}
{\sc Dodig, M.}
\newblock Matrix pencils completions under double rank restrictions.
\newblock {\em Filomat 36}, 4 (2022), 1269--1293.

\bibitem{DoSt09}
{\sc Dodig, M., and Sto$\check{\mbox{s}}$i\'c, M.}
\newblock Similarity class of a matrix with prescribed submatrix.
\newblock {\em Linear Multilinear Algebra 57}, 3 (2009), 217--245.

\bibitem{DoStEJC10}
{\sc Dodig, M., and Sto$\check{\mbox{s}}$i\'c, M.}
\newblock On convexity of polynomial paths and generalized majorizations.
\newblock {\em Electronic Journal of Combinatorics 17}, 1 (2010).

\bibitem{DoSt19}
{\sc Dodig, M., and Sto$\check{\mbox{s}}$i\'c, M.}
\newblock The general matrix completion problem: a minimal case.
\newblock {\em SIAM J. Matrix Anal. Appl. 40}, 1 (2019), 347--369.

\bibitem{FaLe59}
{\sc Farahat, H.~K., and Ledermann, W.}
\newblock Matrices with prescribed characteristic polynomials.
\newblock {\em Proc. Edinburgh Math. Soc. 11}, 3 (1959), 143--146.

\bibitem{Fo75}
{\sc Forney, G.~D.}
\newblock Minimal bases of rational vector spaces with applications to
  multivariable linear systems.
\newblock {\em SIAM Journal Control 13}, 3 (1975), 143--520.

\bibitem{Fri16}
{\sc Friedland, S.}
\newblock {\em {Matrices: algebra, analysis and applications}}.
\newblock World Scientific, Singapore, 2016.

\bibitem{FuSi99}
{\sc Furtado, S., and Silva, F.~C.}
\newblock Embedding a regular subpencil into a general linear pencil.
\newblock {\em Linear Algebra Appl. 295\/} (1999), 61--72.

\bibitem{Gant66}
{\sc Gantmacher, F.}
\newblock {\em {Th\'eorie des matrices, tome 1 et 2}}.
\newblock Dunod, Paris, 1966.

\bibitem{GoKaSch83}
{\sc Gohberg, I., Kaashoek, M.~A., and Schagen, F.~V.}
\newblock Eigenvalues of completion of submatrices.
\newblock {\em Linear Multilinear Algebra 25}, 1 (1983), 55--70.

\bibitem{GLR85}
{\sc Gohberg, I., Lancaster, P., and Rodman, L.}
\newblock {\em Matrix Polynomials}.
\newblock SIAM, Philadelphia, 1985.

\bibitem{HLP88}
{\sc Hardy, G.~H., Littlewood, J.~E., and Polya, G.}
\newblock {\em Inequalities}.
\newblock Cambridge University Press, London, 1988.

\bibitem{HeLMA83}
{\sc Hershkowitz, D.}
\newblock Existence of matrices with prescribed eigenvalues and entries.
\newblock {\em Linear Multilinear Algebra 14\/} (1983), 315--342.

\bibitem{Kail80}
{\sc Kailath, T.}
\newblock {\em Linear Systems}.
\newblock Prentice Hall, New Jersey, 1980.

\bibitem{LoMoZaZaLAA98}
{\sc Loiseau, J., Mondi\'e, S., Zaballa, I., and Zagalak, P.}
\newblock Assigning the {Kronecker} invariants of a matrix pencil by row or
  column completion.
\newblock {\em Linear Algebra Appl. 278\/} (1998), 327--336.

\bibitem{Sa79}
{\sc {Marques de S\`a}, E.}
\newblock Imbedding conditions for $\lambda$-matrices.
\newblock {\em Linear Algebra Appl. 24\/} (1979), 33--50.

\bibitem{Mi58}
{\sc Mirsky, L.}
\newblock Matrices with prescribed characteristic roots and diagonal elements.
\newblock {\em Journal of The London Mathematical Society 33}, 1 (January
  1958), 14--21.

\bibitem{Ro03}
{\sc Roca, A.}
\newblock {\em Asignaci\'on de Invariantes en Sistemas de Control}.
\newblock PhD thesis, Universitat Polit\`ecnica Val\`encia, 2003.

\bibitem{Rose70}
{\sc Rosenbrock, H.~H.}
\newblock {\em State-space and Multivariable Theory}.
\newblock Thomas Nelson and Sons, London, 1970.

\bibitem{Si87}
{\sc Silva, F.~C.}
\newblock Matrices with prescribed eigenvalues and principal submatrices.
\newblock {\em Linear Algebra Appl. 92\/} (1987), 241--250.

\bibitem{Th79}
{\sc Thompson, R.}
\newblock Interlacing inequalities for invariant factors.
\newblock {\em Linear Algebra Appl. 24\/} (1979), 1--31.

\bibitem{TiMe01}
{\sc Tisseur, F., and Meerbergen, K.}
\newblock The quadratic eigenvalue problem.
\newblock {\em SIAM Review 43}, 2 (2001), 235--286.

\bibitem{Vard91}
{\sc Vardulakis, A. I.~G.}
\newblock {\em Linear Multivariable Control}.
\newblock John Wiley and Sons, New York, 1991.

\bibitem{ZaLAA87}
{\sc Zaballa, I.}
\newblock Matrices with prescribed rows and invariant factors.
\newblock {\em Linear Algebra Appl. 87\/} (1987), 113--146.

\bibitem{ZaLAA88}
{\sc Zaballa, I.}
\newblock Interlacing inequalities and control theory.
\newblock {\em Linear Algebra Appl. 101\/} (1988), 9--31.

\bibitem{ZaTi12}
{\sc Zaballa, I., and Tisseur, F.}
\newblock Finite and infinite elementary divisors of matrix polynomials: A
  global approach.
\newblock {\em Manchester Institute for Mathematical Sciences EPrints\/}
  (2012).

\end{thebibliography}
\end{document}